\documentclass[11pt]{article}
\usepackage{amsmath}
\usepackage{amssymb}
\usepackage{amsthm}
\usepackage{latexsym}
\usepackage{hyperref}
\usepackage{enumerate}

\newtheorem{thm}{Theorem}[section]
\newtheorem{cor}[thm]{Corollary}
\newtheorem{lem}[thm]{Lemma}
\newtheorem{prop}[thm]{Proposition}

\newtheorem{thmintro}{Theorem}
\newtheorem{que}[thmintro]{Question}

\providecommand{\norm}[1]{\left\| #1 \right\|}
\newcommand{\enuma}[1]{\begin{enumerate}[\textup{(}a\textup{)}] {#1} \end{enumerate}}

\newcommand{\mh}{\mathbb}
\newcommand{\mr}{\mathrm}
\newcommand{\mc}{\mathcal}
\newcommand{\mf}{\mathfrak}
\newcommand{\ds}{\displaystyle}
\newcommand{\scs}{\scriptstyle}
\newcommand{\ts}{\textstyle}

\newcommand{\N}{\mathbb N}
\newcommand{\Z}{\mathbb Z}
\newcommand{\Q}{\mathbb Q}
\newcommand{\R}{\mathbb R}
\newcommand{\C}{\mathbb C}
\newcommand{\ep}{\epsilon}

\newcommand{\af}{\mr{aff}}

\newcommand{\inp}[2]{\langle #1 \,,\, #2 \rangle}

\begin{document}

\begin{center}
\textbf{\Large Hochschild homology of affine Hecke algebras}\\[3mm]
\Large Maarten Solleveld\\[2mm]
\normalsize Radboud Universiteit Nijmegen \\
Heyendaalseweg 135, 6525AJ Nijmegen, the Netherlands \\
email: m.solleveld@science.ru.nl 
\end{center}

\begin{minipage}{13cm}
\textbf{Abstract.}

Let $\mc H = \mc H (\mc R ,q)$ be an affine Hecke algebra with complex, possibly 
unequal parameters $q$, which are not roots of unity. We compute the Hochschild and the 
cyclic homology of $\mc H$. It turns out that these are independent of $q$ and that they 
admit an easy description in terms of the extended quotient of a torus by a Weyl group, 
both of which are canonically associated to the root datum $\mc R$. 

For positive $q$ we also prove that the representations of the family of algebras 
$\mc H (\mc R ,q^\ep) ,\; \ep \in \C$ come in families which depend analytically on $\ep$. 

Analogous results are obtained for graded Hecke algebras and for Schwartz completions of
affine Hecke algebras.

\texttt{Correction:} in 2021 some problems with the construction of families of 
representations surfaced, these are discussed on pages \pageref{eq:smoothFamily}, 
\pageref{eq:6.2}, \pageref{eq:6.12} and \pageref{eq:6.4}.\\

\end{minipage}

\tableofcontents

\newpage

\section*{Introduction}
\addcontentsline{toc}{section}{Introduction}

The representation theory of affine Hecke algebras has been studied extensively,
for a large part motivated by the connection with reductive $p$-adic groups.
By now this theory is in a very good state, thanks to work of Kazhdan--Lusztig, 
Barbasch--Moy, Delorme--Opdam and many others. Given that the classification of
irreducible representations of affine Hecke algebras is more or less completed
\cite{KaLu,OpSo2,Sol-Irr}, it is natural to look at subtler properties like
extensions, derived categories and homology \cite{BaNi,Nis,OpSo3}.
In this paper we will compute the Hochschild and cyclic homology of affine Hecke 
algebras with possibly unequal parameters. Related results are obtained for graded
Hecke algebras and for Schwartz completions of affine Hecke algebras. This work
can be regarded as a sequel to \cite{SolHomGHA}, where the author determined the
Hochschild homology of graded Hecke algebras (but not the structure as a module
over the centre). 

We discuss the main results in some detail. Let $\mc R = (X,R,Y,R^\vee,\Delta)$ be a 
based root datum with finite Weyl group $W = W(R)$. Let $\mc H = \mc H (\mc R,q)$ be the 
associated affine Hecke algebra over $\C$, with positive (possibly unequal) parameters $q$. 
There is also a Schwartz algebra $\mc S (\mc R,q)$, a topological completion of $\mc H$
which is useful for the study of tempered $\mc H$-representations. We refer to Section 1
for proper definitions of these objects.

Write $T = \mr{Hom}_\Z (X,\C^\times)$ and let $\langle W \rangle$ be a collection of
representatives for the conjugacy classes in $W$. The extended quotient of $T$ by $W$ is
\[
\tilde T / W = \bigsqcup_{w \in \langle W \rangle} T^w / Z_W (w) .
\]
It can be used to parametrize the irreducible $\mc H$-representations  
\cite[Theorem 5.4.2]{Sol-Irr}, in agreement with a conjecture of Aubert, Baum and Plymen
\cite{ABP1}. Then the irreducible $\mc S (\mc R,q)$-representations correspond to
$\tilde T_{un} / W \subset \tilde T / W$, where $T_{un} = \mr{Hom}_\Z (X,S^1)$.

We provide a more precise version of the above parametrization, which at the same time 
extends to the case $q \in \C^\times, |q| \neq 1$. Recall that $\mc H (\mc R,q)$ has a 
commutative subalgebra $\mc A \cong \mc O (T)$, such that $\mc A^W \cong \mc O (T)^W$ 
is the centre of $\mc H (\mc R,q)$. For $w \in \langle W \rangle$ let $T^w_i / Z_w (w) 
,\; 1 \leq i \leq c(w)$ be the connected components of $T^w / Z_W (w)$. 
We consider the family of algebras $\{ \mc H (\mc R ,q^\ep) \mid \ep \in \C \}$.

\begin{thmintro}\label{thmi:1}
(See Theorem \ref{thm:6.6}.) \\
There exist families of $\mc H (\mc R,q^\ep)$-representations 
\[
\{ \pi (w,i,t,\ep) \mid w \in \langle W \rangle, 1 \leq i \leq c(w), t \in T^w_i,
\ep \in \C \}
\]
such that:
\enuma{
\item The representations are irreducible for generic $t$ and $\ep$.
\item The vector space underlying $\pi (w,i,t,\ep)$ depends only on $w$ and $i$. The matrix 
coefficients of this representation are algebraic in $t$ and complex analytic in $\ep$. 
(The latter makes sense because we can identify the vector spaces
$\mc H (\mc R,q^\ep)$ and $\mc H (\mc R,q)$ in a canonical way.)
\item The central character of $\pi (w,i,t,\ep)$ is of the form $W t c_{w,i}^\ep$, where 
$c_{w,i}$ is a homomorphism from $X$ to the subgroup of $\R_{>0}$ generated by the variables
$q_{\alpha^\vee}^{\pm 1/2}$ for all possible coroots $\alpha^\vee$.
\item For $\ep \in \C \setminus \sqrt{-1} \, \R^\times$ the $\mc H (\mc R,q^\ep$)-representations 
$\pi (w,i,t,\ep)$ and $\pi (w',i',t',\ep)$ have the same trace if and only if
$w = w', i = i'$ and $t$ and $t'$ are in the same $Z_W (w)$-orbit.
\item For $\ep \in \C \setminus \sqrt{-1} \, \R^\times$ the collection
\[
\{ \pi (w,i,t,\ep) \mid w \in \langle W \rangle, 1 \leq i \leq c(w), t \in T^w_i / Z_W (w) \}
\]
is a $\Q$-basis of the Grothendieck group of finite dimensional 
$\mc H (\mc R,q^\ep)$-representations.
\item For $\ep \in \R$ the collection
\[
\{ \pi (w,i,t,\ep) \mid w \in \langle W \rangle, 1 \leq i \leq c(w), 
t \in (T^w_i \cap T_{un}) / Z_W (w) \}
\]
is a $\Q$-basis of the Grothendieck group of finite dimensional 
$\mc S (\mc R,q^\ep)$-representations.
}
\end{thmintro}

By analogy with the case of a single parameter $q$ \cite{KaLu}, the author expects that even more
is true. Namely, in (b) the matrix coefficients should depend algebraically on the variables
$q_{\alpha^\vee}^{\pm 1/2}$, while parts (d) and (e) should hold for all but finitely many
$\ep \in \C$. 

The families of representations from Theorem \ref{thmi:1} are our main tool for the computation
of the Hochschild homology. For fixed $w,i,\ep$ the representations $\{ \pi (w,i,t,\ep) \mid
t \in T^w_i \}$ yield an algebra homomorphism
\[
\pi_{w,i,\ep} : \mc H (\mc R ,q^\ep) \to \mc O (T^w_i) \otimes \mr{End}(V_{w,i}) ,
\]
where $V_{w,i}$ is the finite dimensional vector space on which these representations are defined. 
Recall that by the Hochschild--Kostant--Rosenberg theorem 
\[
HH_* \big( \mc O (T^w_i) \otimes \mr{End}(V_{w,i}) \big) \cong \Omega^* (T^w_i) ,
\]
the space of algebraic differential forms on the complex affine variety $T^w_i$.

\begin{thmintro} \label{thmi:2}
(See Theorems \ref{thm:8.2}, \ref{thm:9.1} and \ref{thm:10.1}.) 
\enuma{
\item Let $\ep \in \C \setminus \sqrt{-1} \, \R^\times$. 
The maps $\pi_{w,i,\ep}$ induce an isomorphism
\[
HH_* (\mc H (\mc R,q^\ep)) \to \bigoplus_{w \in \langle W' \rangle} 
\bigoplus_{i=1}^{c(w)} \Omega^* (T^w_i )^{Z_W (w)} \cong
\Omega^* (\tilde T )^W .
\]
The induced action of $Z (\mc H (\mc R,q^\ep)) \cong \mc O (T)^W$ on $\Omega^* (T^w_i)$ is
the same as the action via the embedding
\[
T^w_i \to T : t \mapsto c_{w,i}^\ep t ,
\]
where $c_{w,i}$ is as in Theorem \ref{thmi:1}.c.
\item Let $\ep \in \R$. Part (a) extends to an isomorphism of topological vector spaces
\[
HH_* (\mc S (\mc R,q^\ep)) \to \bigoplus_{w \in \langle W' \rangle} 
\bigoplus_{i=1}^{c(w)} \Omega^*_{sm} (T^w_i \cap T_{un})^{Z_W (w)} \cong
\Omega^*_{sm} (\tilde T_{un} )^W ,
\]
where $\Omega^*_{sm}$ stands for smooth differential forms.
}
\end{thmintro}
Thus the Hochschild homology of $\mc H (\mc R, q^\ep)$ is independent of the parameters
$q^\ep$ and can be expressed in terms of $T$ and $W$. Only the action of the centre 
depends on $q^\ep$, but in a simple way. Another way to look at it is that 
$HH_n (\mc H (\mc R,q^\ep))$ is the space of algebraic $n$-forms on the dual
space of $\mc H (\mc R,q^\ep)$, regarded as a nonseparated variety with a stratification
coming from Theorem \ref{thmi:1}.

All the above results have natural counterparts for graded Hecke algebras. These are
considerably easier to prove and serve as intermediate steps towards Theorems \ref{thmi:1}
and \ref{thmi:2}. One can easily deduce from Theorem \ref{thmi:2} what the (periodic) 
cyclic homology and the Hochschild cohomology of $\mc H (\mc R,q^\ep)$ look like. The
outcome is again that they do not depend on $q^\ep$ and can be expressed with only $T$ and $W$.

There is a second form of Hochschild cohomology,
\[
H^* (\mc H ,\mc H) = \mr{Ext}^*_{\mc H \otimes \mc H^{op}} (\mc H,\mc H) .
\]
It would be quite interesting to compute this, since it is related to deformation
theory and extensions of $\mc H$-bimodules. However, this theory is not dual to
Hochschild homology, so Theorem \ref{thmi:2} says little about it.

In view of Theorems \ref{thmi:1} and \ref{thmi:2} we might wonder how similar two
affine Hecke algebras with the same root datum but different parameter functions are.
From $HH_0$ we see already that $\mc H (\mc R,q)$ and $\mc H (\mc R,q')$ can only be
Morita equivalent if there is an automorphism of $T / W$ that sends every subvariety
$c_{w,i} T^w_i / Z_W (w)$ to a subvariety $c_{w',i'} T^{w'}_{i'} / Z_W (w')$.
It turns out that this condition is rather strong as soon as $R$ contains 
nonperpendicular roots. Indeed, for $\mc R$ of type $\widetilde{A_2}$ it was shown in 
\cite{Yan} that it is only fulfilled if $q' = q$ or $q' = q^{-1}$.

\begin{que} \label{que:1}
Suppose that $\mc H (\mc R,q)$ and $\mc H (\mc R,q')$ are Morita
equivalent or even isomorphic. What are the possibilities for $(q,q')$?
\end{que}

For $\mc R$ of type $\widetilde{A_1}$ it is easily seen that $\mc H (\mc R,q) \cong \C [D_\infty]$
for all $q \in \C \setminus \{-1\}$, where $D_\infty$ is the infinite dihedral group.
But this root datum is exceptional, because it corresponds to the only connected Dynkin diagram
for which the generators of the Coxeter group do not satisfy a braid relation. For irreducible
$\mc R$ of rank at least 2, the aforementioned result of \cite{Yan} suggests that there are only
very few positive answers to Question \ref{que:1}.

Let us briefly describe the organization of the paper. We work in somewhat larger generality than
in the introduction, in the sense that we allow Hecke algebras extended by finite groups of
automorphisms of the root datum. This enhances the applicability, as such algebras appear 
naturally in the representation theory of reductive $p$-adic groups. 

The first section is meant to introduce the notation and to recall several relevant theorems.
In Section 2 we study analytic families of representations of Hecke algebras, in increasing 
generality. This builds upon the author's previous work \cite{Sol-Irr} and leads to
Theorem \ref{thmi:1}. The various homologies of Hecke algebras are computed in Section 3. 
We do it first for graded Hecke algebras, by hand so to say. Subsequently we transfer the
results to affine Hecke algebras and their Schwartz completions, via some arguments 
involving localization at central characters.

We do not study the consequences for the Hochschild homology of reductive $p$-adic groups
here, the author intends to do so in a forthcoming paper. 
\vspace{3mm}

\section{Preliminaries}

\subsection{Affine Hecke algebras}
\label{sec:defAHA}

Let $\mf a$ be a finite dimensional real vector space and let $\mf a^*$ be its dual. 
Let $Y \subset \mf a$ be a lattice and $X = \mr{Hom}_\Z (Y,\Z) \subset \mf a^*$ the dual lattice. Let
\[
\mc R = (X, R, Y ,R^\vee ,\Delta) .
\]
be a based root datum. Thus $R$ is a reduced root system in $X ,\, R^\vee \subset Y$ 
is the dual root system and $\Delta$ is a basis of $R$.
Furthermore we are given a bijection $R \to R^\vee ,\: \alpha \mapsto \alpha^\vee$ such 
that $\inp{\alpha}{\alpha^\vee} = 2$ and such that the corresponding reflections
$s_\alpha : X \to X$ (resp. $s^\vee_\alpha : Y \to Y$) stabilize $R$ (resp. $R^\vee$).
We do not assume that $R$ spans $\mf a^*$. 

The reflections $s_\alpha$ generate the Weyl group $W = W (R)$ of $R$, and 
$S_\Delta := \{ s_\alpha \mid \alpha \in \Delta \}$ is the collection of simple reflections. 
We have the affine Weyl group $W^\af = \mh Z R \rtimes W$ 
and the extended (affine) Weyl group $W^e = X \rtimes W$. Both can be considered as groups
of affine transformations of $\mf a^*$. We denote the translation corresponding to $x \in X$ by $t_x$.
As is well known, $W^\af$ is a Coxeter group, and the basis of $R$ gives rise to a set $S^\af$ 
of simple (affine) reflections. More explicitly, let $\Delta_M^\vee$ be the set of maximal elements of
$R^\vee$, with respect to the dominance ordering coming from $\Delta$. Then
\[
S^\af = S_\Delta \cup \{ t_\alpha s_\alpha \mid \alpha^\vee \in \Delta_M^\vee \} .
\]
We write
\begin{align*}
&X^+ := \{ x \in X \mid \inp{x}{\alpha^\vee} \geq 0
\; \forall \alpha \in \Delta \} , \\
&X^- := \{ x \in X \mid \inp{x}{\alpha^\vee} \leq 0
\; \forall \alpha \in \Delta \} = -X^+ .
\end{align*}
It is easily seen that the centre of $W^e$ is the lattice
\[
Z(W^e) = X^+ \cap X^- .
\]
The length function $\ell$ of the Coxeter system $(W^\af ,S^\af )$ extends naturally to $W^e$.
The elements of length zero form a subgroup $\Omega \subset W^e$ and $W^e = W^\af \rtimes \Omega$.
With $\mc R$ we also associate some other root systems. There is the non-reduced root system
\[
R_{nr} := R \cup \{ 2 \alpha \mid \alpha^\vee \in 2 Y \} .
\]
Obviously we put $(2 \alpha )^\vee = \alpha^\vee / 2$. Let $R_l$
be the reduced root system of long roots in $R_{nr}$:
\[
R_l := \{ \alpha \in R_{nr} \mid \alpha^\vee \not\in 2 Y \} .
\]
We introduce a complex parameter function for $\mc R$ in two equivalent ways.
Firstly, it is a map $q : S^\af \to \mh C^\times$ such that $q(s) = q(s')$ if $s$ and $s'$ 
are conjugate in $W^e$. This extends naturally to a map $q : W^e \to \C^\times$ which is
1 on $\Omega$ and satisfies $q(w w') = q(w) q(w')$ if $\ell (w w') = \ell (w) + \ell (w')$.
Secondly, a parameter function is a $W$-invariant map $q : R_{nr}^\vee \to \C^\times$.
The relation between the two definitions is given by
\begin{equation}\label{eq:parameterEquivalence}
\begin{array}{lll}
q_{\alpha^\vee} = q(s_\alpha) = q (t_\alpha s_\alpha) & \text{if} & \alpha \in R \cap R_l, \\
q_{\alpha^\vee} = q(t_\alpha s_\alpha) & \text{if} & \alpha \in R \setminus R_l, \\
q_{\alpha^\vee / 2} = q(s_\alpha) q(t_\alpha s_\alpha)^{-1} & \text{if} & \alpha \in R \setminus R_l. 
\end{array}
\end{equation}
We speak of equal parameters if $q(s) = q(s') \; \forall s,s' \in S^\af$ and of positive parameters if 
$q(s) \in \R_{>0} \; \forall s \in S^\af$. 

We fix a square root $q^{1/2} : S^\af \to \mh C^\times$.
The affine Hecke algebra $\mc H = \mc H (\mc R ,q)$ is the unique associative
complex algebra with basis $\{ N_w \mid w \in W \}$ and multiplication rules
\begin{equation}\label{eq:multrules}
\begin{array}{lll}
N_w \, N_v = N_{w v} & \mr{if} & \ell (w v) = \ell (w) + \ell (v) \,, \\
\big( N_s - q(s)^{1/2} \big) \big( N_s + q(s)^{-1/2} \big) = 0 & \mr{if} & s \in S^\af .
\end{array}
\end{equation}
In the literature one also finds this algebra defined in terms of the
elements $q(s)^{1/2} N_s$, in which case the multiplication can be described without
square roots. This explains why $q^{1/2}$ does not appear in the notation $\mc H (\mc R ,q)$.

Notice that $N_w \mapsto N_{w^{-1}}$ extends to a $\C$-linear anti-automorphism of $\mc H$,
so $\mc H$ is isomorphic to its opposite algebra.
The span of the $N_w$ with $w \in W$ is a finite dimensional Iwahori--Hecke algebra,
which we denote by $\mc H (W,q)$.

Now we describe the Bernstein presentation of $\mc H$. For $x \in X^+$ we put
$\theta_x := N_{t_x}$. The corresponding semigroup morphism 
$X^+ \to \mc H (\mc R ,q)^\times$ extends to a group homomorphism
\[
X \to \mc H (\mc R ,q)^\times : x \mapsto \theta_x .
\]

\begin{thm}\label{thm:1.1}
\textup{(Bernstein presentation)}
\enuma{
\item The sets $\{ N_w \theta_x \mid w \in W , x \in X \}$ and
$\{ \theta_x N_w \mid w \in W , x \in X \}$ are bases of $\mc H$.
\item The subalgebra $\mc A := \mr{span} \{ \theta_x \mid x \in X \}$
is isomorphic to $\mh C [X]$.
\item The centre of $Z(\mc H (\mc R ,q))$ of $\mc H (\mc R ,q)$ is 
$\mc A^{W}$, where we define the action of $W$ on $\mc A$ by 
$w (\theta_x ) = \theta_{wx}$.
\item For $f \in \mc A$ and $\alpha \in \Delta$
\[
f N_{s_\alpha} - N_{s_\alpha} s_\alpha (f) = 
q(s_\alpha)^{-1/2} \big( f - s_\alpha (f) \big)  (q(s_\alpha) c_\alpha - 1) . 
\]
Here the c-functions are defined as
}
\begin{equation*}
c_\alpha = \left\{
\begin{array}{ll}
{\ds \frac{\theta_\alpha + q(s_\alpha )^{-1/2} q (t_\alpha s_\alpha )^{1/2}}{\theta_\alpha + 1} 
\, \frac{\theta_\alpha - q(s_\alpha )^{-1/2} q (t_\alpha s_\alpha )^{-1/2}}{\theta_\alpha - 1} }
 & \alpha \in R \setminus R_l \\
(\theta_\alpha - q(s_\alpha )^{-1}) (\theta_\alpha - 1)^{-1}
& \alpha \in R \cap R_l .
\end{array}
\right.
\end{equation*}
\end{thm}
\emph{Proof.}
These results are due to Bernstein, see \cite[\S 3]{Lus-Gr}. $\qquad \Box$
\\[2mm]

Consider the complex algebraic torus 
\[
T = \mr{Hom}_{\mh Z} (X, \mh C^\times ) \cong Y \otimes_\Z \C^\times ,
\]
so $\mc A \cong \mc O (T)$ and $Z (\mc H ) = \mc A^{W} \cong \mc O (T / W )$. From Theorem 
\ref{thm:1.1} we see that $\mc H$ is of finite rank over its centre. Let $\mf t = \mr{Lie}(T)$ 
and $\mf t^*$ be the complexifications of $\mf a$ and $\mf a^*$. The direct sum
$\mf t = \mf a \oplus i \mf a$ corresponds to the polar decomposition 
\[
T = T_{rs} \times T_{un} = \mr{Hom}_\Z (X, \R_{>0}) \times \mr{Hom}_\Z (X, S^1)
\]
of $T$ into a real split (or positive) part and a unitary part. The exponential map 
$\exp : \mf t \to T$
is bijective on the real parts, and we denote its inverse by $\log : T_{rs} \to \mf a$.

An automorphism of the Dynkin diagram of the based root system $(R,\Delta )$ is a 
bijection $\gamma : \Delta \to \Delta$ such that
\begin{equation}
\inp{\gamma (\alpha )}{\gamma (\beta )^\vee} = \inp{\alpha}{\beta^\vee} 
\qquad \forall \alpha ,\beta \in \Delta \,.
\end{equation}
Such a $\gamma$ naturally induces automorphisms of $R, R^\vee, W$ and 
$W^\af$. It is easy to classify all diagram automorphisms of $(R ,\Delta)$: they permute the
irreducible components of $R$ of a given type, and the diagram automorphisms of a
connected Dynkin diagram can be seen immediately.

We will assume that the action of $\gamma$ on $W^{\af}$ has been extended in
some way to $W^e$, and then we call it a diagram automorphism of $\mc R$. For example, this 
is the case if $\gamma$ belongs to the Weyl 
group of some larger root system contained in $X$. We regard two diagram automorphisms 
as the same if and only if their actions on $W^e$ coincide.

Let $\Gamma$ be a finite group of diagram automorphisms of $\mc R$ and assume that
$q_{\alpha^\vee} = q_{\gamma (\alpha^\vee)}$ for all $\alpha \in R_{nr}$. Then $\Gamma$
acts on $\mc H$ by algebra automorphisms $\psi_\gamma$ that satisfy
\begin{equation}\label{eq:1.2}
\begin{array}{lll@{\quad}l}
\psi_\gamma (N_w) & = & N_{\gamma (w)} & w \in W , \\
\psi_\gamma (\theta_x) & = & \theta_{\gamma (x)} & x \in X .
\end{array}
\end{equation}
Hence one can form the crossed product algebra $\Gamma \ltimes \mc H = \mc H \rtimes \Gamma$, 
whose natural basis is indexed by the group $(X \rtimes W) \rtimes \Gamma = 
X \rtimes (W \rtimes \Gamma)$. It follows easily from \eqref{eq:1.2} and Theorem 
\ref{thm:1.1}.c that $Z(\mc H \rtimes \Gamma) = {\mc A}^{W \rtimes \Gamma}$. We say that 
the central character of an (irreducible) $\mc H \rtimes \Gamma$-representation is positive 
if it lies in $T^{rs} / (W \rtimes \Gamma)$.
\vspace{2mm}

\subsection{Graded Hecke algebras}

Graded Hecke algebras are also known as degenerate (affine) Hecke algebras.
They were introduced by Lusztig in \cite{Lus-Gr}. We call 
\begin{equation}
\tilde{\mc R} = (\mf a^* ,R, \mf a, R^\vee, \Delta )
\end{equation}
a degenerate root datum. We pick complex numbers $k_\alpha$ for $\alpha \in \Delta$,
such that $k_\alpha = k_\beta$ if $\alpha$ and $\beta$ are in the same $W$-orbit.
The graded Hecke algebra associated to these data is the complex vector space
\[
\mh H = \mh H (\tilde{\mc R},k) = S( \mf t^*) \otimes \C [W] ,
\]
with multiplication defined by the following rules:
\begin{itemize}
\item $\mh C[W]$ and $S (\mf t^* )$ are canonically embedded as subalgebras;
\item for $x \in \mf t^*$ and $s_\alpha \in S$ we have the cross relation
\begin{equation}\label{eq:1.1}
x \cdot s_\alpha - s_\alpha \cdot s_\alpha (x) = 
k_\alpha \inp{x}{\alpha^\vee} .
\end{equation}
\end{itemize}
Multiplication with any $\ep \in \mh C^\times$ defines a bijection $m_\ep : \mf t^* \to \mf t^*$,
which clearly extends to an algebra automorphism of $S(\mf t^* )$. From the cross relation
\eqref{eq:1.1} we see that it extends even further, to an algebra isomorphism
\begin{equation}\label{eq:1.3}
m_\ep : \mh H (\tilde{\mc R},\ep k) \to \mh H (\tilde{\mc R}, k)
\end{equation}
which is the identity on $\mh C[W]$. For $\ep = 0$ this map is well-defined, but
obviously not bijective.

Let $\Gamma$ be a group of diagram automorphisms of $\mc R$, and assume that 
$k_{\gamma (\alpha)} = k_\alpha$ for all $\alpha \in R , \gamma \in \Gamma$.
Then $\Gamma$ acts on $\mh H$ by the algebra automorphisms
\begin{equation}
\begin{split}
& \psi_\gamma : \mh H \to \mh H \,, \\
& \psi_\gamma (x s_\alpha ) = \gamma (x) s_{\gamma (\alpha )} 
  \qquad x \in \mf t^* , \alpha \in \Pi \,.
\end{split}
\end{equation}
By \cite[Proposition 5.1.a]{SolHomGHA} the centre of the resulting crossed product algebra is
\begin{equation}\label{eq:1.4}
Z (\mh H \rtimes \Gamma) = S(\mf t^*)^{W \rtimes \Gamma} = 
\mc O (\mf t / (W \rtimes \Gamma)) .
\end{equation}
We say that the central character of an $\mh H \rtimes \Gamma$-representation is real 
if it lies in $\mf a / (W \rtimes \Gamma)$.
\vspace{2mm}

\subsection{Parabolic subalgebras}

For a set of simple roots $P \subset \Delta$ we introduce the notations
\begin{equation}
\begin{array}{l@{\qquad}l}
R_P = \mh Q P \cap R & R_P^\vee = \mh Q P^\vee \cap R^\vee , \\
\mf a_P = \R P^\vee & \mf a^P = (\mf a^*_P )^\perp ,\\
\mf a^*_P = \R P & \mf a^{P*} = (\mf a_P )^\perp  ,\\
\mf t_P = \C P^\vee & \mf t^P = (\mf t^*_P )^\perp ,\\
\mf t^*_P = \C P & \mf t^{P*} = (\mf t_P )^\perp  ,\\
X_P = X \big/ \big( X \cap (P^\vee )^\perp \big) &
X^P = X / (X \cap \mh Q P ) , \\
Y_P = Y \cap \mh Q P^\vee & Y^P = Y \cap P^\perp , \\
T_P = \mr{Hom}_{\mh Z} (X_P, \mh C^\times ) &
T^P = \mr{Hom}_{\mh Z} (X^P, \mh C^\times ) , \\
\mc R_P = ( X_P ,R_P ,Y_P ,R_P^\vee ,P) & \mc R^P = (X,R_P ,Y,R_P^\vee ,P) , \\
\tilde{\mc R}_P = ( \mf a_P^* ,R_P ,\mf a_P ,R_P^\vee ,P) & 
\tilde{\mc R}^P = (\mf a^*,R_P ,\mf a,R_P^\vee ,P) .
\end{array}
\end{equation}
We define parameter functions $q_P$ and $q^P$ on the root
data $\mc R_P$ and $\mc R^P$, as follows. Restrict $q$ to a function on
$(R_P )_{nr}^\vee$ and determine the value on simple (affine) reflections in 
$W (\mc R_P )$ and $W (\mc R^P )$ by \eqref{eq:parameterEquivalence}. 
Similarly the restriction of $k$ to $P$ is a parameter function for the degenerate 
root data $\tilde{\mc R}_P$ and $\tilde{\mc R}^P$, and we denote it by $k_P$ or $k^P$. 
Now we can define the parabolic subalgebras
\[
\begin{array}{l@{\qquad}l}
\mc H_P = \mc H (\mc R_P ,q_P ) & \mc H^P = \mc H (\mc R^P ,q^P ) , \\
\mh H_P = \mh H (\tilde{\mc R}_P ,k_P ) & \mh H^P = \mh H (\tilde{\mc R}^P ,k^P ) .
\end{array}
\]
We notice that $\mh H^P = S (\mf t^{P*} ) \otimes \mh H_P$, a tensor product of 
algebras. Despite our terminology $\mc H^P$ and $\mc H_P$ are not subalgebras of $\mc H$, 
but they are close. Namely, $\mc H (\mc R^P ,q^P )$ is isomorphic to the subalgebra of 
$\mc H (\mc R ,q)$ generated by $\mc A$ and $\mc H (W (R_P) ,q_P)$. 

We denote the image of $x \in X$ in $X_P$ by $x_P$ and we let $\mc A_P \subset \mc H_P$ 
be the commutative subalgebra spanned by $\{ \theta_{x_P} \mid x_P \in X_P \}$. 
There is natural surjective quotient map
\begin{equation}\label{eq:quotientP}
\mc H^P \to \mc H_P : \theta_x N_w \mapsto \theta_{x_P} N_w .
\end{equation}
Suppose that $\gamma \in \Gamma \ltimes W$ satisfies $\gamma (P) = Q \subseteq \Delta$.
Then there are algebra isomorphisms
\begin{equation}\label{eq:psigamma}
\begin{array}{llcl}
\psi_\gamma : \mc H_P \to \mc H_Q , &
\theta_{x_P} N_w & \mapsto & \theta_{\gamma (x_P)} N_{\gamma w \gamma^{-1}} , \\
\psi_\gamma : \mc H^P \to \mc H^Q , &
\theta_x N_w & \mapsto & \theta_{\gamma x} N_{\gamma w \gamma^{-1}} , \\
\psi_\gamma : \mh H_P \to \mh H_Q  , &
f_P w & \mapsto & (f_P \circ \gamma^{-1}) \gamma w \gamma^{-1}, \\
\psi_\gamma : \mh H_P \to \mh H_Q  , &
f w & \mapsto & (f \circ \gamma^{-1}) \gamma w \gamma^{-1} ,
\end{array}
\end{equation}
where $f_P \in \mc O (\mf t_P)$ and $f \in \mc O (\mf t)$.
Sometimes we will abbreviate $W \rtimes \Gamma$ to $W'$. For example the group
\begin{equation}\label{eq:GammaP}
W'_P := \{ \gamma \in \Gamma \ltimes W \mid \gamma (P) = P \} 
\end{equation}
acts on the algebras $\mc H_P$ and $\mc H^P$. Although $W'_{\Delta} = \Gamma$, for
proper subsets $P \subsetneq \Delta$ the group $W'_P$ need not be contained in $\Gamma$.

To avoid confusion we do not use the notation $W_P$. Instead the parabolic subgroup of $W$ 
generated by $\{ s_\alpha \mid \alpha \in P \}$ will be denoted $W (R_P)$. Suppose that 
$\gamma \in W'$ stabilizes either the root system $R_P$, the lattice $\Z P$ or the vector space
$\Q P \subset \mf a^*$. Then $\gamma (P)$ is a basis of $R_P$, so
$\gamma (P) = w (P)$ and $w^{-1} \gamma \in W'_P$ for a unique $w \in W(R_P)$. Therefore
\begin{equation}\label{eq:WZP}
W'_{\Z P} := \{ \gamma \in W' \mid \gamma (\Z P) = \Z P \} \text{ equals } W(R_P) \rtimes W'_P .
\end{equation}
For $t \in T^P$ and $\lambda \in \mf t^P$ we define an algebra automorphisms
\begin{equation}
\begin{array}{llll@{\quad}lll@{\quad}l}
\phi_t : & \mc H^P & \to & \mc H^P, & \phi_t (\theta_x N_w) & = & 
t (x) \theta_x N_w & x \in X, w \in W, \\
\phi_\lambda : & \mh H^P & \to & \mh H^P, & \phi_\lambda (f h) & = & f(\lambda) f h & 
f \in S (\mf t^{P*}), h \in \mh H_P .
\end{array}
\end{equation}
For $t \in K_P := T^P \cap T_P$ this descends to an algebra automorphism 
\begin{equation}\label{eq:twistKP}
\psi_t : \mc H_P \to \mc H_P , \quad \theta_{x_P} N_w \mapsto t(x_P) \theta_{x_P} N_w
\qquad t \in K_P .
\end{equation}
We can regard any representation $(\sigma ,V_\sigma)$ of $\mc H (\mc R_P ,q_P )$ as a 
representation of $\mc H^P = \mc H (\mc R^P ,q^P)$ via the quotient map \eqref{eq:quotientP}. 
Thus we can construct the $\mc H$-representation
\[
\pi (P,\sigma ,t) := \mr{Ind}_{\mc H^P}^{\mc H} (\sigma \circ \phi_t ) .
\]
Representations of this form are said to be parabolically induced.
Similarly, for any $\mh H_P$-representation $(\rho, V_\rho)$ and any $\lambda \in \mf t^P$
there is an $\mh H^P$-representation $\rho \circ \phi_\lambda$.
The corresponding parabolically induced representation is
\[
\pi (P,\rho,\lambda) := \mr{Ind}_{\mh H^P}^{\mh H} (\rho \circ \phi_\lambda) .
\]
In case we include a group of diagram automorphisms $\Gamma$, we will also use the representations
\begin{align*}
& \pi^\Gamma (P,\sigma ,t) := 
\mr{Ind}_{\mc H^P}^{\mc H \rtimes \Gamma} (\sigma \circ \phi_t ) , \\
& \pi^\Gamma (P,\rho,\lambda) := 
\mr{Ind}_{\mh H^P}^{\mh H \rtimes \Gamma} (\rho \circ \phi_\lambda) .
\end{align*}

\subsection{Lusztig's reduction theorems} 
\label{sec:redthm}

The study of irreducible representations of $\mc H \rtimes \Gamma$  is simplified by 
two reduction theorems, which are essentially due to Lusztig \cite{Lus-Gr}.
The first one reduces to the case of modules whose central character is positive on
the lattice $\Z R_l$. The second one relates these to modules of an associated graded
Hecke algebra. 

Given $t \in T$ and $\alpha \in R$, \cite[Lemma 3.15]{Lus-Gr} tells us that
\begin{equation}\label{eq:sAlphaFixt}
s_\alpha (t) = t \text{ if and only if } \alpha (t) = 
\begin{cases} 1 & \text{if } \alpha^\vee \notin 2 Y \\
\pm 1 & \text{if } \alpha^\vee \in 2 Y .
\end{cases}
\end{equation}
We define $R_t := \{ \alpha \in R \mid s_\alpha (t) = t \}$. The collection of long
roots in $R_{t,nr}$ is $\{ \beta \in R_l \mid \beta (t) = 1\}$. Let $F_t$ be the unique
basis of $R_t$ consisting of roots that are positive with respect to $\Delta$.
We can define a parameter function $q_t$ for the based root datum
\[
\mc R_t := (X,R_t,Y,R_t^\vee,F_t)
\]
via restriction from $R_{nr}^\vee$ to $R_{t,nr}^\vee$. Furthermore we write
\[
P(t) := \Delta \cap \Q R_t .
\]
Then $R_{P(t)}$ is a parabolic root subsystem of $R$ that contains $R_t$ 
as a subsystem of full rank. Let $t = u c \in T_{un} T_{rs}$ be the polar
decomposition of $t \in T$. We note that $R_{uc} \subset R_u$, that 
$W'_{uc} \subset W'_u$ and that the lattice
\[
\Z P(t) = \Z R \cap \Q R_u 
\]
can be strictly larger than $\Z R_t$. We will phrase the first reduction theorem such that 
it depends mainly on the unitary part $u$ of $t$, it will decompose a representation in parts
corresponding to the point of the orbit $W' u$.

For a finite set $U \subset T / W$, let $Z_U (\mc H) \subset Z (\mc H)$ be the ideal of functions 
vanishing at $U$. Let $\widehat{Z (\mc H )_U}$ be the formal completion of $Z (\mc H)$ with
respect to the powers of the ideal $Z_U (\mc H)$ and define
\begin{equation}\label{eq:ZHt}
\widehat{\mc H_U} = \widehat{Z (\mc H )_U} \otimes_{Z (\mc H)} \mc H .
\end{equation}
Similarly, for $t \in T$ let $\widehat{\mc A_t}$ denote the formal completion of $\mc A \cong 
\mc O (T)$ with respect to the powers of the ideal $\{ f \in \mc O (T) \mid f(t) = 0 \}$.
Inside $\widehat{\mc H_U}$ we have the formal completion of $\mc A$ at $U$, which 
by the Chinese remainder theorem is isomorphic to
\begin{equation}
\widehat{\mc A_U} := \widehat{Z (\mc H )_U} \otimes_{\mc A^W} \mc A \cong
\bigoplus_{t \in U} \widehat{\mc A_t} .
\end{equation}
In this notation we can rewrite
\[
\widehat{Z (\mc H )_{W t}} \cong \widehat{\mc A_{W t}}^{W} \cong \widehat{\mc A_t}^{W_t} .
\]
Analogous statements hold for $\mc H \rtimes \Gamma$.
Given a subset $\varpi \subset W' t$ we let $1_\varpi \in \widehat{\mc A_{W' t}}$
be the idempotent corresponding to $\bigoplus_{s \in \varpi} \widehat{\mc A_s}$.

\begin{thm}\label{thm:2.1}
\textup{(First reduction theorem)} \\
There is a natural isomorphism of $\widehat{Z (\mc H \rtimes \Gamma)_{W' uc}}$-algebras
\[
\widehat{\mc H (\mc R^{P(u)}, q^{P(u)})_{W'_{\Z P (u)} uc}} \rtimes W'_{P(u)} \cong 
1_{W'_{\Z P(u)} uc} \big( \widehat{\mc H_{W' uc}} \rtimes \Gamma \big) 1_{W'_{\Z P(u)} uc}
\]
It can be extended (not naturally) to isomorphism of 
$\widehat{Z (\mc H \rtimes \Gamma)_{W' uc}}$-algebras
\[
\widehat{\mc H_{W' uc}} \rtimes \Gamma  \cong  
M_{[W' : W'_{\Z P(u)}]} \Big( 1_{W'_{\Z P(u)} uc} \big( \widehat{\mc H_{W' uc}} \rtimes \Gamma \big)  
1_{W'_{\Z P(u)} uc} \Big) ,
\]
where $M_n (A)$ denotes the algebra of $n \times n$-matrices with coefficients in an algebra~$A$. 
In particular the algebras
\[
\widehat{\mc H (\mc R^{P(u)}, q^{P(u)})_{W'_{\Z P(u)} uc}} \rtimes W'_{P(u)} \quad \text{and} \quad
\widehat{\mc H_{W' uc}} \rtimes \Gamma 
\]
are Morita equivalent.
\end{thm}
\emph{Proof.}
This is a variation on \cite[Theorem 8.6]{Lus-Gr}. Compared to Lusztig we substituted his $R_{uc}$ 
by a larger root system, we replaced the subgroup $Y \otimes \langle v_0 \rangle \subset T$ by
$T_{rs} = Y \otimes \R_{>0}$ and we included the automorphism group~$\Gamma$. These changes are 
justified in the proof of \cite[Theorem 2.1.2]{Sol-Irr}. $\quad \Box$
\\[2mm]

By \eqref{eq:sAlphaFixt} we have $\alpha (u) = 1$ for all $\alpha \in R_l \cap \Q R_u$, so
$\alpha (t) = \alpha (u) \alpha (c) > 0$ for such roots. Hence Theorem \ref{thm:2.1} allows us 
to restrict our attention to $\mc H \rtimes \Gamma$-modules whose central character is positive
on the sublattice $\Z R_l \subseteq X$.

By definition $u$ is fixed by $W'_u \supset W (R_{P(u)})$, so the map
\begin{equation}\label{eq:expMap}
\exp_u : \mf t \to T ,\;  \lambda \mapsto u \exp (\lambda)
\end{equation}
is $W'_u$-equivariant. 

Analogous to \eqref{eq:ZHt}, let $Z_{W \lambda} (\mh H) \subset Z (\mh H)$ be the maximal ideal
of functions vanishing at $W \lambda \in \mf t / W$. Let $\widehat{Z (\mh H )_{W \lambda}}$ 
be the formal completion of $Z (\mh H)$ with respect to $Z_{W \lambda} (\mh H)$ and define 
\begin{equation}
\widehat{\mh H_{W \lambda}} := \widehat{Z (\mh H )_{W' \lambda}} \otimes_{Z (\mh H)} \mh H .
\end{equation}
The corresponding formal completion of $S (\mf t^*) \cong \mc O (\mf t)$ is
\[
\widehat{S (\mf t^*)_{W \lambda}} := \widehat{Z (\mh H )_{W \lambda}} 
\otimes_{S (\mf t^*)^W} S (\mf t^*) \cong \bigoplus_{\mu \in W \lambda} \widehat{S (\mf t^*)_\mu} .
\]
The map \eqref{eq:expMap} induces a $W'_u$-equivariant isomorphism
\[
\widehat{\mc A_{u \exp \lambda}} \to \widehat{S (\mf t^*)_{\lambda}} : f \mapsto f \circ \exp_u ,
\]
which restricts to an isomorphism
\begin{equation}\label{eq:Phiu}
\Phi_{u,W' \lambda} : \widehat{Z (\mc H \rtimes \Gamma)_{W' u \exp (\lambda)}} 
\to \widehat{Z (\mh H \rtimes \Gamma)_{W' \lambda}} .
\end{equation}
We define a parameter function $k_u$ for the degenerate root datum $\tilde{\mc R}_u$ by
\begin{equation}\label{eq:2.3}
2 k_{u,\alpha} = \log q_{\alpha^\vee / 2} + (1 + \alpha (u)) \log q_{\alpha^\vee}
\qquad \alpha \in R_u .
\end{equation}
Let $q^{\Z / 2}$ be the subgroup of $\mh C^\times$ generated by
$\{ q_{\alpha^\vee}^{\pm 1/2} \mid \alpha^\vee \in R_{nr}^\vee \}$. 

\begin{thm}\label{thm:2.2}
\textup{(Second reduction theorem)} \\
Suppose that 1 is the only root of unity in $q^{\Z / 2}$.
Let $u \in T_{un}^{W'}$ and let $\lambda \in \mf t$ be such that
\begin{equation}\label{eq:2.2}
\inp{\alpha}{\lambda}, \inp{\alpha}{\lambda} + k_{u,\alpha} 
\notin \pi i \mh Z \setminus \{0\} \qquad \forall \alpha \in R \,.
\end{equation} 
\enuma{
\item The map \eqref{eq:Phiu} extends to an algebra isomorphism 
\[
\Phi_{u,W' \lambda} : \widehat{\mc H_{W' u \exp (\lambda)}} \rtimes \Gamma \to 
\widehat{\mh H (\tilde{\mc R},k_u)_{W' \lambda}} \rtimes \Gamma.
\]
\item The algebras $\widehat{\mc H_{W' u \exp (\lambda)}} \rtimes \Gamma$ and
$\widehat{\mh H (\tilde{\mc R}^{P(u)}, k^{P(u)})_{W'_{\Z P(u)} \lambda}} \rtimes W'_{P(u)}$
are Morita equivalent.
}
\end{thm}
\emph{Proof.}
For part (a) see \cite[Theorem 9.3]{Lus-Gr}. Our conditions on $q$ replace the 
assumption \cite[9.1]{Lus-Gr}. Part (b) follows from part (a) and Theorem \ref{thm:2.1}. 
$\qquad \Box$ \vspace{2mm}

\subsection{Schwartz algebras}
\label{sec:temp}

An important tool to study $\mc H$-representations is restriction to the commutative 
subalgebra $\mc A \cong \mc O (T)$. We say that $t \in T$ is a weight of $(\pi,V)$ if there
exists a $v \in V \setminus \{ 0 \}$ such that $\pi (a) v = a(t) v$ for all $a \in \mc A$.
Temperedness of $\mc H$-representations, which is defined via $\mc A$-weights, is analogous to 
temperedness of representations of reductive groups. Via the Langlands classification for 
affine Hecke algebras \cite[Section 2.2]{Sol-Irr} these representations are essential 
in the classification of irreducible representations.

The antidual of $\mf a^{*+} :=  \{ x \in \mf a^* \mid  \inp{x}{\alpha^\vee} \geq 0 
\: \forall \alpha \in \Delta \}$ is
\begin{equation} 
\mf a^- = \{ \lambda \in \mf a \mid \inp{x}{\lambda} \leq 0 \: \forall x \in \mf a^{*+} \} = 
\big\{ \sum\nolimits_{\alpha \in \Delta} 
\lambda_\alpha \alpha^\vee \mid \lambda_\alpha \leq 0 \big\} .
\end{equation}
The interior $\mf a^{--}$ of $\mf a^-$ equals $\big\{ {\ts \sum_{\alpha \in \Delta}} 
\lambda_\alpha \alpha^\vee \mid \lambda_\alpha < 0 \big\}$
if $\Delta$ spans $\mf a^*$, and is empty otherwise. We write $T^- = \exp (\mf a^-)$ and
$T^{--} = \exp (\mf a^{--})$. 

Let $t = |t| \cdot t |t|^{-1} \in T_{rs} \times T_{un}$ be the polar decomposition of $t$.
A finite dimensional $\mc H$-representation is called tempered if $|t| \in T^-$ for all its 
$\mc A$-weights $t$, and anti-tempered if $|t|^{-1} \in T^-$ for all such $t$. 

We say that an irreducible  $\mc H$-representation belongs to the discrete 
series (or simply: is discrete series) if all its $\mc A$-weights lie in $T^{--} T_{un}$. 
In particular the discrete series is empty if $\Delta$ does not span $\mf a^*$.
The discrete series is a starting point for the construction of irreducible representations:
they can all be realized as a subrepresentation of the induction of a discrete series 
representation of a parabolic subalgebra of $\mc H$.
The notions tempered and discrete series apply equally well to $\mc H \rtimes \Gamma$,
since that algebra contains $\mc A$ and the action of $\Gamma$ on $T$ preserves $T^-$.

This terminology also extends naturally to graded Hecke algebras, via the $S (\mf t^*)$-weights
of a representation. Thus we say that a finite dimensional $\mh H \rtimes \Gamma$-representation
is tempered if all its $S (\mf t^*)$-weights lie in $\mf a^- + i \mf a$ and we say that it
is discrete series if it is irreducible and all its $S (\mf t^*)$-weights lie in 
$\mf a^{--} + i \mf a$. By construction Lusztig's two reduction theorems from Section 
\ref{sec:redthm} preserve the properties temperedness and discrete series.

\begin{prop}\label{prop:5.1}
Let $P \subset \Delta$.
\enuma{
\item Let $\sigma$ be a finite dimensional $\mc H_P \rtimes \Gamma_P$-representation. 
For $t \in T^P$ the $\mc H \rtimes \Gamma$-representation
$\mr{Ind}_{\mc H_P \rtimes \Gamma_P}^{\mc H \rtimes \Gamma} (\sigma \circ \phi_t)$
is tempered if and only if $t \in T^P_{un}$ and $\sigma$ is tempered.
\item Let $\rho$ be a finite dimensional $\mh H_P \rtimes \Gamma_P$-representation. 
For $\lambda \in \mf t^P$ the $\mh H \rtimes \Gamma$-representation
$\mr{Ind}_{\mh H_P \rtimes \Gamma_P}^{\mh H \rtimes \Gamma} (\rho \circ \phi_\lambda)$
is tempered if and only if $\lambda \in i \mf a^P$ and $\rho$ is tempered.
}
\end{prop}
\emph{Proof.}
For (a) see \cite[Lemma 3.1.1.b]{Sol-Irr} and for (b) see \cite[Lemma 2.2]{SolGHA}. 
$\qquad \Box$ \\[2mm]

Now we will recall the construction of the Schwartz algebra $\mc S$ of an affine Hecke 
algebra \cite{DeOp1}, which we will use in Section \ref{sec:homS}.
For this we assume that the parameter function $q$ is positive, because otherwise
there would not be a good link with $C^*$-algebras. As a topological vector space $\mc S$ 
will consist of rapidly decreasing functions on $W^e$, with respect to a suitable length 
function $\mc N$. For example we can take a $W$-invariant norm on $X \otimes_\Z \R$ and put
$\mc N (w t_x) = \norm{x}$ for $w \in W$ and $x \in X$. Then we can define, for $n \in \N$, 
the following norm on $\mc H$:
\[
p_n \big( \sum_{w \in W^e} h_w N_w \big) = \sup_{w \in W^e} |h_w| (\mc N (w) + 1)^n .
\]
The completion of $\mc H$ with respect to these norms is the Schwartz algebra 
$\mc S = \mc S (\mc R,q)$. It is known from \cite[Section 6.2]{Opd-Sp} that it is a 
Fr\'echet algebra. Diagram automorphisms of $\mc R$ induce automorphisms of $\mc S$, so the
crossed product algebra $\mc S \rtimes \Gamma$ is well-defined. By \cite[Lemma 2.20]{Opd-Sp} 
a finite dimensional $\mc H \rtimes \Gamma$-representation is 
tempered if and only if it extends continuously to an $\mc S \rtimes \Gamma$-representation. 

Next we describe the Fourier transform for $\mc H$ and $\mc S$.
An induction datum for $\mc H$ is a triple $(P,\delta,t)$, with $P \subset \Delta , t \in T^P$ 
and $\delta$ a discrete series representation of $\mc H_P$. A large part of the representation 
theory of $\mc H$ is built on the parabolically induced representations
\[
\pi (P,\delta,t) = \mr{Ind}_{\mc H^P}^{\mc H} (\delta \circ \phi_t)  \quad \text{and} \quad 
\pi^\Gamma (P,\delta,t) = \mr{Ind}_{\mc H^P}^{\mc H \rtimes \Gamma} (\delta \circ \phi_t) .
\]
Let $\Xi$ be the space of all such induction data $(P,\delta,t)$, with $\delta$ up to
equivalence of $\mc H_P$-representations. It carries a natural structure of a complex affine 
variety with finitely many components of different dimensions. Furthermore $\Xi$ has 
a compact submanifold
\[
\Xi_{un} := \big\{ (P,\delta,t) \in \Xi \mid t \in T^P_{un} \big\} ,
\]
which is a disjoint union of finitely many compact tori.

Given a discrete series representation $(\delta ,V_\delta)$ of a parabolic subalgebra 
$\mc H_P$ of $\mc H$, the $\mc H \rtimes \Gamma$-representation $\pi^\Gamma (P,\delta,t)$
can be realized on the vector space $V_\delta^\Gamma := \C [\Gamma W^P] \otimes V_\delta$,
which does not depend on $t \in T^P$. Here $W^P$ is a specified set of representatives for
$W / W(R_P)$. Let $\mc V_\Xi^\Gamma$ be the vector bundle over $\Xi$ whose fiber at 
$\xi = (P,\delta,t)$ is $V_\delta^\Gamma$, and let 
\[
\mc O (\Xi; \mr{End}(\mc V_\Xi^\Gamma)) := \bigoplus\nolimits_{P,\delta} \mc O (T^P) \otimes
\mr{End}_\C (V_\delta^\Gamma)
\]
be the algebra of polynomials sections of the endomorphism bundle $\mr{End} (\mc V_\Xi^\Gamma)$.
The Fourier transform for $\mc H \rtimes \Gamma$ is
\begin{align*}
& \mc F : \mc H \rtimes \Gamma \to \mc O (\Xi; \mr{End}(\mc V_\Xi^\Gamma)), \\
& \mc F (h) (\xi) = \pi^\Gamma (\xi) (h) .
\end{align*}
It extends to an algebra homomorphism
\begin{equation}
\mc F : \mc S \rtimes \Gamma \to C^\infty (\Xi_{un}; \mr{End}(\mc V_\Xi^\Gamma )) :=
\bigoplus\nolimits_{P,\delta} C^\infty (T^P_{un}) \otimes \mr{End}_\C (V_\delta^\Gamma) ,
\end{equation}
defined by the same formula. 

We will need a groupoid $\mc G$ over the power set of $\Delta$, defined as follows. 
For $P,Q \subset \Delta$ the collection of arrows from $P$ to $Q$ is
\[
\mc G_{PQ} = \{ (g,u) \in \Gamma \ltimes W \times K_P \mid g(P) = Q \} .
\]
Whenever it is defined, the multiplication in $\mc G$ is 
\[
(g',u') \cdot (g,u) = (g' g, g^{-1}(u) u') .
\]
This groupoid acts from the left on $\Xi$ by
\begin{equation}\label{eq:GactsXi}
(g,u) (P,\delta,t) = \big( g(P),\delta \circ \psi^{-1}_u \circ \psi^{-1}_g, g(ut) \big) .
\end{equation} 
This is the projection of an action of $\mc G$ on parabolically induced representations via
intertwining operators. These operators provide an action of $\mc G$ on $C^\infty (\Xi_{un}; 
\mr{End}(\mc V_\Xi^\Gamma ))$, see \cite[Theorem 4.33]{Opd-Sp} and \cite[(3.16)]{Sol-Irr}. 
The Plancherel isomorphism for (extended) affine Hecke algebras with positive parameters reads:

\begin{thm} \label{thm:5.2}
The Fourier transform $\mc F$ induces an isomorphism of Fr\'echet algebras
\[
\mc S \rtimes \Gamma \to C^\infty (\Xi_{un}; \mr{End}(\mc V_\Xi^\Gamma ))^{\mc G} .
\]
\end{thm}
\emph{Proof.}
See \cite[Theorem 5.3]{DeOp1} and \cite[Theorem 3.2.2]{Sol-Irr}. $\qquad \Box$
\vspace{2mm}

From this isomorphism we see in particular that there are unique central idempotents
$e_{P,\delta} \in \mc S \rtimes \Gamma$ such that 
\begin{equation} \label{eq:5.1}
\begin{aligned}
& e_{P,\delta} \mc S \rtimes \Gamma \cong 
\big( C^\infty (T^P_{un}) \otimes \mr{End}_\C (V_\delta^\Gamma) \big)^{\mc G_{P,\delta}} ,\\
& \mc G_{P,\delta} := \big\{ g \in \mc G \mid g (P,\delta,T^P_{un}) = (P,\delta,T^P_{un}) \big\} .
\end{aligned}
\end{equation}
For a suitable collection $\mc P$ of pairs $(P,\delta)$ we obtain decompositions
\begin{equation} \label{eq:5.2}
\begin{aligned}
& \mc S \rtimes \Gamma && = && 
\bigoplus\nolimits_{(P,\delta) \in \mc P} e_{P,\delta} \mc S \rtimes \Gamma ,\\
& Z (\mc S \rtimes \Gamma) && = && 
\bigoplus\nolimits_{(P,\delta) \in \mc P}  C^\infty (T^P_{un})^{\mc G_{P,\delta}} .
\end{aligned}
\end{equation}
Via \eqref{eq:5.1} the subalgebra $e_{P,\delta} \mc H \rtimes \Gamma$ of 
$\mc S \rtimes \Gamma$ is isomorphic to a subalgebra of
$\big( \mc O (T^P) \otimes \mr{End}_\C (V_\delta^\Gamma) \big)^{\mc G_{P,\delta}}$.
We note that $e_{P,\delta} Z (\mc H \rtimes \Gamma )$ is isomorphic to the restriction of
$Z (\mc H \rtimes \Gamma) \cong \mc O (T)^{W'}$ to the image of $(P,\delta,T^P_{un})$
in $T / W'$. We will see from Lemma \ref{lem:6.3} that it is just $\mc O (T^P)^{\mc G_{P,\delta}}$.

Clearly $e_{P,\delta} \mc H \rtimes \Gamma$ is dense in $e_{P,\delta} \mc S \rtimes \Gamma$,
so for any closed ideal $I \subset e_{P,\delta} \mc S \rtimes \Gamma$ of finite codimension
the canonical maps
\begin{equation} \label{eq:5.3}
e_{P,\delta} \mc H \rtimes \Gamma / \big( I \cap e_{P,\delta} \mc H \rtimes \Gamma \big) \to
\big( e_{P,\delta} \mc H \rtimes \Gamma + I \big) / I \to
e_{P,\delta} \mc S \rtimes \Gamma / I 
\end{equation}
are isomorphisms of $e_{P,\delta} \mc H \rtimes \Gamma$-modules. 

\begin{lem}\label{lem:5.3}
The multiplication map
\[
\mu_{P,\delta} : e_{P,\delta} Z (\mc S \rtimes \Gamma ) \otimes_{e_{P,\delta} Z (\mc H 
\rtimes \Gamma )} e_{P,\delta} \mc H \rtimes \Gamma \to e_{P,\delta} \mc S \rtimes \Gamma 
\]
is an isomorphism of $e_{P,\delta} Z (\mc S \rtimes \Gamma )$-modules.
\end{lem}
\emph{Proof.}
Clearly $\mu_{P,\delta}$ is $e_{P,\delta} Z (\mc S \rtimes \Gamma )$-linear. By Theorem
\ref{thm:5.2} $e_{P,\delta} \mc S \rtimes \Gamma$ is of finite rank over 
$e_{P,\delta} Z (\mc S \rtimes \Gamma )$ and we know that $e_{P,\delta} \mc H \rtimes \Gamma$ 
is dense in $e_{P,\delta} \mc S \rtimes \Gamma$. Thus the image of $\mu_{P,\delta}$ is
closed and dense in $e_{P,\delta} \mc S \rtimes \Gamma$, which means that $\mu_{P,\delta}$ 
is surjective. 

Let $I_\xi$ be the maximal ideal of $e_{P,\delta} Z (\mc S \rtimes \Gamma ) \cong
C^\infty (T^P_{un})^{\mc G_{P,\delta}}$ corresponding to $\xi = (P,\delta,t)$. Applying 
\eqref{eq:5.3} to $I_\xi^n e_{P,\delta} \mc S \rtimes \Gamma$ we see that
$\mu_{P,\delta}$ induces isomorphisms
\begin{equation}\label{eq:5.4}
e_{P,\delta} Z (\mc S \rtimes \Gamma ) / I_\xi^n \otimes_{e_{P,\delta} Z (\mc H 
\rtimes \Gamma )} e_{P,\delta} \mc H \rtimes \Gamma \to 
e_{P,\delta} \mc S \rtimes \Gamma / I_\xi^n (e_{P,\delta} \mc S \rtimes \Gamma ) ,
\end{equation}
for all $n \in \Z_{>0}$. Consider any $x \in \ker \mu_{P,\delta}$. Its annihilator Ann$(x)$
is a closed ideal of $e_{P,\delta} Z (\mc S \rtimes \Gamma )$ and $e_{P,\delta} Z(\mc S \rtimes
\Gamma) x \cong e_{P,\delta} Z(\mc S \rtimes \Gamma) / \mr{Ann}(x)$. By \eqref{eq:5.4}
$I_\xi^n / I_\xi^n \mr{Ann}(x) = 0$ for all $\xi \in (P,\delta,T^P_{un})$. Therefore
Ann$(x)$ is not contained in any of the closed maximal ideals $I_\xi$. Consequently
Ann$(x) = e_{P,\delta} Z (\mc S \rtimes \Gamma )$ and $x = 0. \qquad \Box$
\vspace{3mm}

\section{Analytic families of representations}

\subsection{Positive parameters}
\label{sec:positive}

In this subsection we assume that $q$ is positive. We will define analytic families of 
representations and show that they can be used to construct the dual space of $\mc H (\mc R,q)$.

\begin{thm} \label{thm:6.1}
Let $\xi, \xi' \in \Xi$.
\enuma{
\item The $\mc H \rtimes \Gamma$-representations $\pi^\Gamma (\xi)$ and $\pi^\Gamma (\xi')$
have the same trace if and only if there exists a $g \in \mc G$ such that $g \xi = \xi'$.
\item Suppose that $\xi, \xi' \in \Xi_{un}$. Then $\pi^\Gamma (\xi)$ and $\pi^\Gamma (\xi')$ are 
completely reducible, and they have a common irreducible subquotient if only if there exists a 
$g \in \mc G$ such that $g \xi = \xi'$. Moreover $\pi^\Gamma (\xi) \cong \pi^\Gamma (\xi')$ 
in this situation.
}
\end{thm}
\emph{Proof.}
(a) \cite[Lemma 3.1.7]{Sol-Irr} provides the "if"-part. 
By \cite[page 44]{Sol-Irr} the "only if"-part can be reduced to certain positive induction data, 
to which \cite[Theorem 3.3.1]{Sol-Irr} applies.\\
(b) See Corollary 3.1.3 and Theorem 3.3.1 of \cite{Sol-Irr}. We note that this is essentially a 
consequence of Theorem \ref{thm:5.2}. $\qquad \Box$ 
\\[2mm]

Since $\gamma (P) \subset \Delta$ for all $\gamma \in \Gamma$ and $P \subset \Delta$,
\eqref{eq:GactsXi} describes an action of $\Gamma$ on $\Xi$. 
Given an induction datum $\xi = (P,\delta,t) \in \Xi$ we put
\[
\Gamma_\xi = \{ \gamma \in \Gamma_P \mid \delta \circ \phi_t \cong \delta \circ \phi_t
\circ \psi_\gamma^{-1} \text{ as } \mc H^P \text{-representations }\} .
\]
Now we fix $(P,\delta ,u) \in \Xi_{un}$ and we let $\sigma$ be an irreducible direct summand of
$\mr{Ind}_{\mc H^P}^{\mc H^P \rtimes \Gamma_{(P,\delta,u)}} (\delta \circ \phi_u)$.
In this case we abbreviate $\Gamma_{(P,\delta,u)} = \Gamma_\sigma$. By Clifford theory
the representations 
$\mr{Ind}_{\mc H^P}^{\mc H^P \rtimes \Gamma_\sigma} (\delta \circ \phi_u)$ and 
$\mr{Ind}_{\mc H^P \rtimes \Gamma_\sigma}^{\mc H^P \rtimes \Gamma_P} (\sigma)$
are completely reducible.

Let $T^\sigma$ be the connected component of $T^{W (R_P) \rtimes \Gamma_\sigma}$
that contains $1 \in T$. Notice that $T^\sigma \subset T^P$ because $T_P^{W(R_P)}$ is finite.
We call
\begin{equation}\label{eq:smoothFamily}
\big\{ \pi_\sigma (t) = \mr{Ind}_{\mc H^P \rtimes \Gamma_\sigma}^{\mc H \rtimes \Gamma}
(\sigma \circ \phi_t) \mid t \in T^\sigma \big\} 
\end{equation}
an analytic $d$-dimensional family of $\mc H \rtimes \Gamma$-representations, where $d$ is the 
dimension of the complex algebraic variety $T^\sigma$. By Proposition \ref{prop:5.1} the 
representations $\pi_\sigma (t)$ are tempered if and only if 
$t \in T^\sigma_{un} := T^\sigma \cap T_{un}$. We refer to 
$\{ \pi_\sigma (t) \mid t \in T^\sigma_{un} \}$ as a tempered analytic family
or an analytic family of $\mc S$-representations. \\[1mm]

\texttt{Correction (2021):} To guarantee that the union of analytic families of representations 
spans $G_\Q (\mc H \times \Gamma)$, we need to allow more general $\sigma$. Namely, we need all 
irreducible tempered elliptic representations of a parabolic subalgebra $\mc H^P \rtimes \Gamma'_P$, 
where $\Gamma'_P$ is an arbitrary subgroup of $\Gamma_P$. Here elliptic means that the image in 
$G_\Q (\mc H^P \rtimes \Gamma'_P)$ does not belong to the span of the representations induced from 
proper parabolic subalgebras of $\mc H^P \rtimes \Gamma'_P$.
Such a $\sigma$ is a direct summand of the $\mc H^P \rtimes \Gamma'_P$-representation associated to a
triple $(P,\delta,u) \in \Xi_{un}$. With that, the below arguments hold in the necessary generality.\\

Since there are only finitely many pairs $(P,\delta)$ and since two $u$'s in the same 
$T^\sigma$-coset give rise to the same analytic family, there exist only finitely many analytic 
families of $\mc H \rtimes \Gamma$-representations. Recall from \cite[Theorem 2.58]{OpSo2}
that the absolute value of any $\mc A_P$-weight of $\delta$ is a monomial in the variables 
$q (s_\alpha)^{\pm 1/2}$. That is, it lies in 
\begin{equation}
W (R_P) u_\delta Y_P \otimes_\Z q^{\Z / 2}
\end{equation}
where $u_\delta \in T_{P,un}$ is the unitary part of such a weight.
Hence all $\mc A$-weights of $\pi_\sigma (t)$ lie in
\begin{equation}\label{eq:monomial}
W' t u  u_\delta Y \otimes_\Z q^{\Z / 2} . 
\end{equation}
We note that the $\mc H^P \rtimes \Gamma_\sigma$-representation $\pi_\sigma (t)$ is only
defined for $t \in T^{W(R_P) \rtimes \Gamma_\sigma}$, so it is impossible to extend 
\eqref{eq:smoothFamily} to a larger connected subset of $T^P$. As discussed in \cite{Sol-Irr} after
(3.30), the representations $\pi_\sigma (t)$ are irreducible for $t$ in a Zariski-open dense
subset of $T^\sigma$. Like in \eqref{eq:5.1} the group
\begin{equation} \label{eq:Gsigma}
\mc G_\sigma := \big\{ g \in \mc G \mid g(P) = P, \delta \circ \psi_g^{-1} \cong \delta, 
g \text{ stabilizes } \{ \pi_\sigma (t) \mid t \in T^\sigma_{un} \} \big\}
\end{equation}
acts on $T^\sigma$. (Here the advantage of considering only $T^\sigma_{un}$ is that Theorem 
\ref{thm:6.1}.b applies.) By Theorem \ref{thm:6.1}.a $\pi_\sigma (t)$ and $\pi_\sigma (g t)$ 
have the same trace for all $t \in T^{\sigma}$ and $g \in \mc G_\sigma$. We will see later 
that every representation $\pi_{\sigma}(t')$ with $t'$ in a different $\mc G_\sigma$-orbit 
has a different trace.

Let $G (A)$ denote the Grothendieck group of finite dimensional $A$-representations, for suitable
algebras or groups $A$. As explained in \cite[Section 3.4]{Sol-Irr}, 
\begin{equation*}
G_\Q (\mc H \rtimes \Gamma) := \Q \otimes_\Z G (\mc H \rtimes \Gamma) 
\end{equation*}
can be built from analytic families of representations (called smooth families in that paper). 
By this we mean that we can choose a collection of analytic families 
$\{ \pi_{\sigma_i} (t) \mid t \in T^{\sigma_i} \}$ such that the set
\begin{equation}\label{eq:6.1}
\bigcup\nolimits_i  \{ \pi_{\sigma_i} (t) \mid t \in T^{\sigma_i}  / \mc G_{\sigma_i} \}
\end{equation}
spans $G_\Q (\mc H \rtimes \Gamma)$. In \cite[Section 3.4]{Sol-Irr} this is actually done with
"Langlands constituents" of $\pi_{\sigma_i} (t)$. But by \cite[Lemma 2.2.6]{Sol-Irr} the
Langlands formalism is such that the other constituents of $\pi_{\sigma_i}(t)$ are smaller in
a certain sense. Hence the non-Langlands constituents are already accounted for by families of 
smaller dimension.

All this can be made much more concrete for the algebra 
\[
\mc H (\mc R ,1) \rtimes \Gamma = \C [X \rtimes (W \rtimes \Gamma)] = \mc O (T) \rtimes W' .
\] 
By classical results which go back to Frobenius and Clifford, its irreducible representations 
with $\mc O (T)^{W'}$-character $W' t$ are in bijection with the irreducible representations 
$\sigma$ of the isotropy group $W'_t$, via
\begin{equation}\label{eq:6.9}
\sigma \mapsto \mr{Ind}_{X \rtimes W'_t}^{X \rtimes W'} (\C_t \otimes \sigma) .
\end{equation}
With this one can easily determine the dual space of $X \rtimes W'$. Let
\[
\tilde T = \{ (w,t) \in T \times W' \mid w(t) = t \} .
\]
Then $W'$ acts on $\tilde T$ by $w \cdot (w',t) = (w w' w^{-1}, w(t))$ and $\tilde T / W'$ 
is called the extended quotient of $T$ by $W'$. Let $\langle W' \rangle$ be a set of 
representatives for the conjugacy classes in $W'$. Then
\begin{equation}\label{eq:6.10}
\tilde T / W' \cong \bigsqcup_{w \in \langle W' \rangle} T^w / Z_{W'}(w) ,
\end{equation}
where $Z_{W'}(w)$ denotes the centralizer of $w$ in $W'$. Let $T_i^w / Z_{W'}(w) ,
i=1 ,\ldots, c(w)$ be the connected components of $T^w / Z_{W'}(w)$, so that
\begin{equation}\label{eq:6.11}
\tilde T / W' \cong \bigsqcup_{w \in \langle W' \rangle} 
\bigsqcup_{1 \leq i \leq c(w)} T_i^w / Z_{W'}(w) .
\end{equation}
By the above there exists a continuous bijection from $\tilde T / W'$ to the space of irreducible 
complex representations of $X \rtimes W'$. We vary a little and take the 
Grothendieck group of $X \rtimes W'$ instead. Let $C_w \subset W'$ be the cyclic group generated 
by $w$ and let $\rho_w$ be a $C_w$-representation.
This yields an analytic family of $X \rtimes W'$-representations 
\[
\big\{ \pi_{w,i} (t) := \mr{Ind}_{X \rtimes C_w}^{X \rtimes W'} 
(\C_t \otimes \rho_w) \mid t \in T_i^w \big\} .
\]
For generic $t$, $\pi_{w,i}(t)$ is irreducible if $W'_t = C_w$. Moreover $\pi_{w,i}(t)$ arise
by induction from an elliptic representation of $X \rtimes (W(R_P) \rtimes \Gamma'_P)$, where
$W(R_P) \rtimes \Gamma'_P$ is a parabolic subgroup of $W \rtimes \Gamma$ minimally containing $w$.

By Artin's Theorem \cite[Theorem 17]{Ser} it is possible to choose the $\rho_w$ such that $G_\Q (W'_t)$ 
is spanned by $\{ \mr{Ind}_{C_w}^{W'_t}(\rho_w) \mid w \in W'_t \}$, \\[1mm]

\texttt{Correction (2021): The next claim is false, it can fail for\\
instance if $W' = \Z / 4 \Z$.} \\

and a basis is obtained by 
taking only one $w$ from every conjugacy class in $W'_t$. It follows that the representations
\begin{equation} \label{eq:6.2}
\big\{ \pi_{w,i} (t) \mid w \in \langle W' \rangle, 
1 \leq i \leq c(w), t \in T_i^w / Z_{W'}(w) \big\}
\end{equation}
form a basis of $G_\Q (X \rtimes W')$. Notice that the space underlying \eqref{eq:6.2} is exactly
the extended quotient $\tilde T / W'$, since all the involved representations are different. 
Having described this Grothendieck group conveniently, we recall how it can be compared with the 
Grothendieck groups of associated affine Hecke algebras.

\begin{thm} \label{thm:6.2}
\textup{\cite[Section 2.3]{Sol-Irr}} \\
Let $q$ be positive. There exists a natural $\Q$-linear bijection 
\[
\zeta^\vee : G_\Q (\mc H \rtimes \Gamma) \to G_\Q (X \rtimes W')
\]
such that:
\begin{enumerate}[\textup{(}1\textup{)}]
\item $\zeta^\vee$ restricts to a bijection between the corresponding Grothendieck groups of 
tempered representations;
\item $\zeta^\vee$ commutes with parabolic induction, in the sense that for a subgroup 
$\Gamma'_P \subset \Gamma_P$, a tempered $\mc H^P \rtimes \Gamma'_P$-representation 
$\pi$ and $t \in T^{W(R_P) \rtimes \Gamma'_P}$ we have
\[
\zeta^\vee \big( \mr{Ind}_{\mc H^P \rtimes \Gamma'_P}^{\mc H \rtimes \Gamma} 
(\pi \circ \phi_t) \big) = \mr{Ind}_{X \rtimes (W (R_P) \rtimes \Gamma'_P 
)}^{X \rtimes (W \rtimes \Gamma)} (\zeta^\vee (\pi) \circ \phi_t) ;
\]
\item if $u \in T_{un}$ and $\pi$ is a virtual representation with central character in 
$W' u T_{rs}$, then so is $\zeta^\vee (\pi)$.
\end{enumerate}
\end{thm}

Let $d \in \N$. By property (2) $\zeta^\vee$ sends analytic $d$-dimensional families of 
$\mc H \rtimes \Gamma$-representations to analytic $d$-dimensional families of 
$X \rtimes W'$-representations - except that the latter need not be generically irreducible. 
The image can not be part of a higher
dimensional analytic family, for that would violate the bijectivity of $\zeta^\vee$. 
Since no $d$-dimensional variety is a finite union of varieties of dimension strictly smaller 
than $d$, and since there are only finitely many analytic families, $\zeta^\vee$ restricts to 
a bijection between the subspaces (of the Grothendieck groups) spanned by analytic families of 
dimension at least $d$. Let us call these subspaces $G_\Q^d (\mc H \rtimes \Gamma)$ and 
$G_\Q^d (X \rtimes W')$. Then $\zeta^\vee$ induces a bijection 
\begin{equation}\label{eq:6.12}
G_\Q^d (\mc H \rtimes \Gamma) / G_\Q^{d+1} (\mc H \rtimes \Gamma) \to
G_\Q^d (X \rtimes W') / G_\Q^{d+1} (X \rtimes W') .
\end{equation}
By \eqref{eq:6.2} the right-hand side has a basis which is parametrized by the $d$-dimensional 
part of $\widetilde T / W'$. \\[1mm]

\texttt{Correction (2022):} Unfortunately the below construction of families of representations
does not work in general, because of the problems with \eqref{eq:6.2} and because an application
of Theorem \ref{thm:6.2} could result in virtual $\mc H \rtimes \Gamma$-representations.
Nevertheless a collection of families of representations satisfying Lemma \ref{lem:6.3} below
can be found in most cases:
\begin{itemize}
\item When the parameters $q$ are of ``geometric type", one can use families of standard 
representations from \cite[Paragraph 3.3]{AMS3}.
\item For most positive-valued $q$ one can use the technique with parameter deformations from
\cite[proof of Lemma 6.4]{SolHecke}, to reduce to the previous case.
\item We are not aware of any examples of graded Hecke algebras for which it is clear that they 
do not possess families of representations satisfying Lemma \ref{lem:6.3}.
\end{itemize}
In all these cases, the remainder of the paper is valid without further changes.\\[2mm]

We can use this to pick a suitable collection of analytic families of 
$\mc H \rtimes \Gamma$-representations. Fix a $d$-dimensional component of $\widetilde T / W'$
and let $C$ be its image in $T / W'$. In the notation of \eqref{eq:6.11} let $\tilde C$ 
be the collection of $(w,i)$ such that $w \in \langle W' \rangle ,\; 1 \leq i \leq c(w)$ and 
$T^w_i / W' = C$. By Theorem \ref{thm:6.2} and \eqref{eq:6.12} there exist analytic families of
$\mc H \rtimes \Gamma$-representations $\{ \pi_{\sigma_{w,i}}(t) \mid (w,i) \in \tilde C , 
t \in T^{\sigma_{w,i}} \}$ such that
\[
\zeta^\vee \big( G_\Q^{d+1} (\mc H \rtimes \Gamma) \big) + \zeta^\vee \big( \{ \pi_{\sigma_{w,i}}(t) 
\mid (w,i) \in \tilde C , t \in T^{\sigma_{w,i}} \} \big)
\]
contains a basis of 
\[
\text{span} \{ \pi \in \text{Irr}(X \rtimes W') \mid c \text{ is the central character of } \pi \} 
\subset G_\Q (X \rtimes W') ,
\]
for all generic points $c \in C$. This procedure yields analytic $d$-dimensional families whose 
span in $G_\Q (\mc H \rtimes \Gamma)$ is mapped bijectively to $G_\Q^d (X \rtimes W') \big/ 
G_\Q^{d+1} (X \rtimes W')$ by $\zeta^\vee$. By \eqref{eq:Gsigma} we may take the parameters for 
$\pi_{\sigma_{w,i}}$ in $T^{\sigma_{w,i}} / \mc G_{\sigma_{w,i}}$, when we consider these 
representations in $G (\mc H \rtimes \Gamma)$. Let us list some properties of
\begin{equation} \label{eq:6.3}
\{ \pi_{\sigma_{w,i}}(t) \mid w \in \langle W' \rangle , 1 \leq i \leq c(w) ,
t \in T^{\sigma_{w,i}} / \mc G_{\sigma_{w,i}} \} .
\end{equation}

\begin{lem}\label{lem:6.3}
\enuma{
\item The parameter space $T^{\sigma_{w,i}} / \mc G_{\sigma_{w,i}}$ is isomorphic to 
$T^w_i / Z_{W'}(w)$, as orbifolds.
\item $\pi_{\sigma_{w,i}}(t)$ and $\pi_{\sigma_{w,i}}(t')$ have the same trace if and only if
$t$ and $t'$ lie in the same $\mc G_{\sigma_{w,i}}$-orbit.
\item The family of representations \eqref{eq:6.3} contains no repetitions 
and forms a basis of $G_\Q (\mc H \rtimes \Gamma)$.
}
\end{lem}
\noindent \emph{Proof.}
By the above construction and since different $w \in W'$ act differently on $T$,
\begin{align*}
& \zeta^\vee \big( G_\Q^{d+1} (\mc H \rtimes \Gamma) \big) + 
\zeta^\vee \big( \text{span} \{ \pi_{\sigma_{w,i}}(t) \mid t \in T^{\sigma_{w,i}} \} \big) = \\
& \zeta^\vee \big( G_\Q^{d+1} (X \rtimes W') \big) + 
\zeta^\vee \big( \text{span} \{ \pi_w (t) \mid t \in T^w_i / Z_{W'}(w) \} \big) .
\end{align*}
Hence the parameter space $T^w_i / Z_{W'}(w)$ is diffeomorphic to $T^{\sigma_{w,i}}$
modulo the equivalence relation 
\[
t \sim t' \Longleftrightarrow \pi_{\sigma_{w,i}}(t) = 
\pi_{\sigma_{w,i}}(t') \text{ in } G_\Q (\mc H \rtimes \Gamma) .
\]
By Theorem \ref{thm:6.1}.b $t \sim t'$ if and only if there exists a $g \in \mc G$ with 
$g (P,\delta,ut) = (P,\delta,ut')$, where $\pi_{\sigma_{w,i}}$ is a direct summand of 
$\pi (P,\delta,u)$. By definition any $g \in \mc G_{\sigma_{w,i}}$ yields such equivalences 
for all $t \in T^{\sigma_{w,i}}$, so $T^{\sigma_{w,i}} / \sim$ is a quotient of 
$T^{\sigma_{w,i}} / \mc G_{\sigma_{w,i}}$. On the other hand, elements 
$h \in \mc G_{(P,\delta)} \setminus \mc G_{\sigma_{w,i}}$ can only produce such
equivalences on lower dimensional subsets of $T^{\sigma_{w,i}}$. But any such relation which 
does not already come from $\mc G_{\sigma_{w,i}}$ would destroy the orbifold structure of 
$T^{\sigma_{w,i}} / \mc G_{\sigma_{w,i}}$, so 
\[
T^w_i / Z_{W'}(w) \cong T^{\sigma_{w,i}} / \sim \, \cong T^{\sigma_{w,i}} / \mc G_{\sigma_{w,i}} .
\]
This proves (a) and since we work in $G_\Q (\mc H \rtimes \Gamma)$, it also implies (b).

We noted that \eqref{eq:6.2} contains no repetitions and is a basis of $G_\Q (X \rtimes W')$. 
Since the collection of analytic families 
\[
\{ \zeta^\vee (\pi_{\sigma_{w,i}}(t)) \mid t \in w \in \langle W' \rangle , 
1 \leq i \leq c(w) , t \in T^{\sigma_{w,i}} / \mc G_{\sigma_{w,i}} \}
\]
is parametrized by the same space $\tilde T / W'$, it also enjoys these properties. 
As $\zeta^\vee$ is bijective, these two properties already hold for \eqref{eq:6.3} in 
$G_\Q (\mc H \rtimes \Gamma) . \qquad \Box$
\vspace{2mm}

\subsection{Complex parameters}
\label{sec:complex}

Everything we discussed in the previous subsection has a natural analogue for graded 
Hecke algebras. The condition $q$ positive translates to $k$ real valued, but with the 
scaling isomorphisms \eqref{eq:1.3} we can also reach complex parameters.

In the graded setting we have to replace $T^P$ und $T^P_{un}$ by the vector spaces 
$\mf t^P$ and $i \mf a^P$. In the groupoid $\mc G$ we have to omit the parts $K_P$, 
so we get a groupoid $\mc W'$ with
\[
\mc W'_{PQ} = \{ g \in \Gamma \ltimes W \mid g(P) = Q \} \qquad P,Q \subset \Delta .
\]

\begin{thm}\label{thm:6.4}
Suppose that $k_\alpha \in \R$ for all $\alpha \in R$. 
There exists a natural $\Q$-linear bijection
\[
\zeta^\vee : G_\Q (\mh H \rtimes \Gamma) \to G_\Q (S (\mf t^*) \rtimes W')
\]
such that:
\begin{enumerate}[\textup{(}1\textup{)}]
\item $\zeta^\vee (\pi)$ is a tempered virtual $S (\mf t^*) \rtimes W'$-representation 
if and only if $\pi$ is a tempered virtual $\mh H \rtimes \Gamma$-representation.
\item $\zeta^\vee$ commutes with parabolic induction, in the sense that for a subgroup 
$\Gamma'_P \subset \Gamma_P$, a tempered $\mh H^P \rtimes \Gamma'_P$-representation 
$\pi$ and $\lambda \in \mf t^{W(R_P) \rtimes \Gamma'_P}$ we have
\[
\zeta^\vee \big( \mr{Ind}_{\mh H^P \rtimes \Gamma'_P}^{\mh H \rtimes \Gamma} 
(\pi \circ \phi_\lambda) \big) = 
\mr{Ind}_{S (\mf t^*) \rtimes (W (R_P) \rtimes \Gamma'_P)}^{S (\mf t^*) 
\rtimes (W \rtimes \Gamma)} (\zeta^\vee (\pi) \circ \phi_\lambda) ;
\]
\item if $\lambda \in i \mf a$ and $\pi$ is a virtual representation with central character in 
$W' \lambda + \mf a$, then so is $\zeta^\vee (\pi)$;
\item if $\pi$ is tempered and admits a real central character, then 
$\zeta^\vee (\pi) = \pi \circ m_0$ with $m_0$ as in \eqref{eq:1.3}.
\end{enumerate}
\end{thm}
\emph{Proof.}
For $\mh H \rtimes \Gamma$-representations which admit a central character whose component
in $i \mf a$ is small, this follows from Theorems \ref{thm:6.2} and \ref{thm:2.2}. The extra
condition (4) is actually the starting point for the construction of $\zeta^\vee$, see 
\cite[Section 2.3]{Sol-Irr}. We can generalize this to all $\mh H \rtimes \Gamma$-representations
via the scaling isomorphisms
\[
m_\ep : \mh H (\tilde{\mc R},\ep k) \rtimes \Gamma \to \mh H (\tilde{\mc R},k) \rtimes \Gamma
\]
for $\ep > 0 . \qquad \Box$
\\[2mm]

An induction datum for $\mh H$ is a triple $(P,\delta,\lambda)$ such that 
$P \subset \Delta, \delta$ is a discrete series representation of $\mh H_P$ and 
$\lambda \in \mf t^P$. Suppose that $\mu \in i \mf a^P$ and that $\rho$ is an 
irreducible direct summand of the representation 
$\mr{Ind}_{\mh H^P}^{\mh H^P \rtimes \Gamma_{(P,\delta,\mu)}}(\delta \circ \phi_\mu)$ (which
is completely reducible). We abbreviate $\Gamma_\rho = \Gamma_{(P,\delta,\mu)}$
and $\mf t^\rho = (\mf t^P)^{\Gamma_\rho}$. We say that 
\begin{equation} \label{eq:6.13}
\big\{ \pi_\rho (\lambda) := \mr{Ind}_{\mh H^P \rtimes \Gamma_\rho}^{\mh H \rtimes \Gamma}
(\rho \circ \phi_\lambda) \mid \lambda \in \mf t^\rho \big\}
\end{equation}
is an analytic family of $\mh H \rtimes \Gamma$-representations. Notice that 
\begin{equation}\label{eq:6.6}
\pm \mu \in \mf t^\rho 
\end{equation}
and that the central character of $\pi_\rho (-\mu)$ equals the central character of 
$\pi^\Gamma (P,\delta,0) = \mr{Ind}_{\mh H^P}^{\mh H \rtimes \Gamma} \delta$. Moreover the
weights satisfy an integrality property with respect to 
\[
\Z k := \Z \{ k_\alpha \mid \alpha \in R \} . 
\]
By \cite[Theorem A.7]{Opd-Sp} every $S (\mf t_P^*)$-weight $\nu$ of $\delta$ is a residual 
point for $(\tilde{\mc R}_P, k_P)$. The classification of such distinguished points 
in \cite[Section 4]{HeOp} shows that $\alpha (\nu) \in \Z k_P$ for all $\alpha \in R_P$.
By \cite[Theorem 6.4]{BaMo2} this property is preserved under induction, that is, 
$\alpha (\nu) \in \Z k$ for every $\alpha \in R$ and every $S (\mf t^*)$-weight $\nu$ of 
$\mr{Ind}_{\mh H^P}^{\mh H \rtimes \Gamma} \delta$. We conclude that all 
$S (\mf t^*)$-weights of $\pi_\rho (\lambda)$ lie in
\begin{equation}\label{eq:6.4}
W' (\lambda + \mu) + \mr{Hom}(\Z R, \Z k) ,
\end{equation}
where we consider $\mr{Hom}(\Z R, \Z k)$ as a lattice in $\mf a_R \subset \mf t$.
 
To the analytic family $\pi_\rho$ we associate the group
\[
\mc W'_\rho = \big\{ g \in \mc W' \mid g(P) = P , \delta \circ \psi_g^{-1} \cong \delta , 
g \text{ stabilizes } \{ \pi_\rho (\lambda) \mid \lambda \in i \mf a^\rho \} \big\}. 
\]
\texttt{Correction (2022):} The new comments on pages \pageref{eq:smoothFamily}, 
\pageref{eq:6.2} and \pageref{eq:6.12} 
also apply to the below construction and to Corollary \ref{cor:6.5}. If we can instead use one 
of the good cases listed on page \pageref{eq:6.12}, the remainder of the paper holds for those
families of $\mh H \rtimes \Gamma$-representations. Otherwise, there is still a satisfactory 
substitute for Corollary \ref{cor:6.5} in terms of families of virtual $\mh H \rtimes 
\Gamma$-representations, worked out in \cite{SolHHtwisted}.\\[2mm]

\noindent
Our discussion of $\C [X] \rtimes W' = \mc O (T) \rtimes W'$ before Theorem \ref{thm:6.2} 
applies just as well to $\mh H(\tilde{\mc R},0) \rtimes \Gamma = S (\mf t^*) \rtimes W'$.
The considerations that lead to Lemma \ref{lem:6.3} also show:

\begin{cor}\label{cor:6.5}
There exist analytic families of $\mh H \rtimes \Gamma$-representations $\pi_{\rho_w}$ for
$w \in \langle W' \rangle$, such that:
\enuma{
\item The parameter space $\mf t^{\rho_w} / \mc W'_{\rho_w}$ is isomorphic to
$\mf t^w / Z_{W'}(w)$, as orbifolds.
\item $\pi_{\rho_w} (\lambda)$ and $\pi_{\rho_w}(\lambda')$ have the same trace if and only
if $\lambda$ and $\lambda'$ lie in the same $\mc W'_{\rho_w}$-orbit.
\item The family of representations 
\[
\{ \pi_{\rho_w}(\lambda ) \mid w \in \langle W' \rangle, 
\lambda \in \mf t^{\rho_w} / \mc W'_{\rho_w} \}
\]
contains no repetitions and forms a basis of $G_\Q (\mh H \rtimes \Gamma)$.
}
\end{cor}

The scaling isomorphisms $m_\ep$ for $\ep \in \C^\times$ from \eqref{eq:1.3} provide an
immediate generalization of Theorem \ref{thm:6.4} and Corollary \ref{cor:6.5} to graded
Hecke algebras $\mh H (\tilde R,\ep k) \rtimes \Gamma$ such that $k$ is real valued and
$\ep \in \C$. To retain the formulation of Theorem \ref{thm:6.4} one must modify the 
notions tempered and discrete series, in the sense that everywhere in Section \ref{sec:temp}
one must use $\ep \mf a^-$ instead of $\mf a^-$. Let us call this $\ep$-tempered and
$\ep$-discrete series. We can also consider the family of algebras 
\[
\big\{ \mh H (\tilde{\mc R} , \ep k) \rtimes \Gamma \mid \ep \in \C \big\} .
\]
For every $\ep \in \C$ the analytic families of 
$\mh H (\tilde{\mc R} , \ep k) \rtimes \Gamma$-representations 
\[
\{ \pi_{m_\ep^* (\rho_w )} (\lambda) \mid w \in \langle W' \rangle, 
\lambda \in \mf t^{\rho_w} / \mc W'_{\rho_w} \}
\]
fulfill Corollary \ref{cor:6.5}. We note that the matrix coefficients of
$\pi_{m_\ep^* (\rho_w )} (\lambda)$ depend algebraically on $\ep$ and $\lambda$. More 
specifically, by \eqref{eq:6.4} the $S(\mf t^*)$-weights of $\pi_{m_\ep^* (\rho_w)} (-\ep \mu)$ 
lie in $\mr{Hom} (\Z R, \Z k \ep) \subset \mf t_R$ and depend linearly on $\ep$.

Next we generalize Section \ref{sec:positive} to certain
complex parameter functions $q$. As the analogue of a line $\ep \R$ we use a one-parameter
subgroup of $\C^\times$. We consider $q^\ep$ for a positive parameter function $q$ and 
$\ep \in \C$. Then \eqref{eq:2.3} becomes
\[
k_{u,\ep,\alpha} = \big( \log q_{\alpha^\vee / 2} + (1 + \alpha (u)) \log q_{\alpha^\vee} \big)
\ep / 2 \qquad \alpha \in R_u .
\]
Notice that $\ep = 0$ corresponds to the algebra $\mc H (\mc R ,q^0) \rtimes \Gamma = 
\C [X] \rtimes W'$, which we understand very well. On the other hand, for $\ep \in i \R 
\setminus \{0\}$ the parameter $q_{\alpha^\vee}$ can be a root of unity (unequal to 1), 
in which case the affine Hecke algebra $\mc H (\mc R ,q)$ can differ substiantially from 
those with positive parameters. Therefore we assume that $\Re (\ep) \neq 0$.

\begin{thm}\label{thm:6.6}
Theorem \ref{thm:6.2} and Lemma \ref{lem:6.3} also hold for parameter functions $q^\ep$ with
$q$ positive and $\ep \in \C \setminus i \R$.
For the interpretation of Theorem \ref{thm:6.2} we have to replace $\mf a^-$ by $\ep \mf a^-$,
$T_{rs}$ by $\exp (\ep \mf a)$ and tempered by $\ep$-tempered. 

Moreover we can choose the required analytic families of $\mc H (\mc R,q^\ep) \rtimes 
\Gamma$-representations such that the matrix coefficients of the representations depend 
analytically on $\ep$ and all $\mc A$-weights depend polynomially on the variables
$\big\{ q_{\alpha^\vee}^{\pm \ep/2} \mid \alpha^\vee \in R_{nr}^\vee \big\}$.
\end{thm}
\emph{Proof.}
Every $t \in T$ can be written uniquely as $t = u c$ with $u \in T_{un}$ and 
$c \in \exp (\ep \mf a)$, since $\Re (\ep) \neq 0$. In this setting Theorem \ref{thm:2.1} 
remains valid. Moreover $k_{u,\ep,,\alpha} \in \ep \R$ and 
$\langle \alpha, \lambda \rangle \in \ep \R$ for all $\alpha \in R$ and $\lambda \in \ep \mf a$, 
so $\langle \alpha, \ep \mf a \rangle$ and $\langle \alpha, \ep \mf a \rangle + k_{u,\ep,\alpha}$ 
have positive distance to $\pi i \Z \setminus \{0\}$. This ensures that Theorem \ref{thm:2.2} 
applies to all $\lambda$ in a suitable tubular neighborhood of $\ep \mf a$ in $\mf t$. 
Its part (b) tells us that
$\widehat{\mc H (\mc R,q^\ep)}_{W'u \exp (\lambda)} \rtimes \Gamma$ is Morita equivalent to 
$\widehat{\mh H (\tilde{\mc R}^{P(u)}, \ep k^{P(u)})}_{W'_{\Z P(u)} \lambda} \rtimes W'_{P(u)}$,
while $\mh H (\tilde{\mc R}^{P(u)}, \ep k^{P(u)})$ is isomorphic to a graded Hecke algebra
with real parameter function $k^{P(u)}$ by means of $m_{\ep}$. To 
$\mh H (\tilde{\mc R}^{P(u)}, k^{P(u)})$ we can apply Theorem \ref{thm:6.4}, which gives us
Theorem \ref{thm:6.2} for $\mc H \rtimes \Gamma$-representations with central character
in $W' u \exp (\ep \mf a)$. Combining these results for all $u \in T_{un}$ yields the map
\[
\zeta^\vee : G (\mc H (\mc R ,q^\ep ) \rtimes \Gamma) \to G (X \rtimes W') .
\]
It makes Theorem \ref{thm:6.2} true for $q^\ep$, because it does so on all sets 
$u \exp (\ep \mf a)$.

Let us construct the required analytic families of $\mc H (\mc R,q^\ep ) \rtimes 
\Gamma$-representations. We start with \eqref{eq:6.3} and let $\sigma$ be any of the
$\sigma_{w,i}$. With Theorem \ref{thm:2.2} we transfer the $\mc H^P \rtimes \Gamma_\sigma
$-representation to a tempered representation $\tilde \sigma$ of $\mh H^{P(u)} \rtimes 
\tilde{\Gamma}$, for a suitable subgroup $\tilde \Gamma \subset W'$. Taking for $u$ the 
unitary part of an $\mc A$-weight of $\sigma$, we can achieve that $\tilde \sigma$ has 
real central character. Then $m_\ep^* (\tilde \sigma)$ is an $\ep$-tempered representation of 
$\mh H (\tilde{\mc R}^P(u), \ep k^{P(u)}) \rtimes \tilde{\Gamma}$. Now Theorem \ref{thm:2.2}
in the opposite direction yields an $\ep$-tempered representation of 
$\mc H (\mc R^P , (q^P)^\ep ) \rtimes \Gamma_\sigma$, which we call $m_\ep^* (\sigma)$.
We claim that the analytic families of $\mc H (\mc R ,q^\ep) \rtimes \Gamma$-representations
\[
\{ \pi_{m_\ep^* (\sigma_{w,i})}(t) \mid w \in \langle W' \rangle , 1 \leq i \leq c(w) ,
t \in T^{\sigma_{w,i}} / \mc G_{\sigma_{w,i}} \} 
\]
fulfill Lemma \ref{lem:6.3}. By construction $\zeta^\vee (\pi_{m_\ep^* (\sigma_{w,i})}(t))$ 
equals $\zeta^\vee (\sigma_{w,i})(t)$. In view of Lemma \ref{lem:6.3} and
since Theorem \ref{thm:6.2} holds for $q$ and $q^\ep$, this implies the claim.

In view of \cite[(4.6)]{Sol-Irr}, $m_\ep^* (\sigma)$ is also the composition of $\sigma$ 
with the algebra homomorphism $\rho_\ep$ from \cite[Proposition 4.1.2]{Sol-Irr}. Although 
in \cite{Sol-Irr} the author assumed that $\ep \in [-1,1]$, the relevant calculations 
remain valid for $\ep \in \C \setminus i \R$. This shows that the matrix coefficients of 
$m_\ep^* (\sigma)$ depend analytically on $\ep$. More precisely, for every $v \in W^e$ the 
operator $\pi_{m_\ep^* (\sigma_{w,i})}(t,N_v)$ depends in a complex analytic way on 
$\ep$ and $t$.

By \eqref{eq:monomial} all $\mc A$-weights of $\pi_{\sigma_{w,i}}(t)$ are of the form
$\gamma (t u c)$ where $\gamma \in W'$ and $uc$ is a weight of $\sigma_{w,i}$ with unitary 
part $u \in T_{un}$ and absolute value $c \in Y \otimes_\Z q^{\Z / 2}$. This in turn 
implies that all $\mc A$-weights of $\pi_{m_\ep^* (\sigma_{w,i})}(t)$ are of the form
\begin{equation} \label{eq:6.weights}
\gamma (t u c^\ep )  \qquad c^\ep \in q^{\ep \Z / 2},
\end{equation}
and in particular that (for fixed $t$) they depend polynomially on
$\big\{ q_{\alpha^\vee}^{\pm \ep/2} \mid \alpha^\vee \in R_{nr}^\vee \big\}. \qquad \Box$ 
\\[2mm]
Since the whole situation depends in complex analytic fashion on $\ep$, it is natural to
expect that Theorem \ref{thm:6.6} extends to all but finitely many $\ep \in \C$. 
Indeed this is known when $q$ is an equal parameter function \cite{KaLu}, but the methods 
of this paper are insufficient to prove it.
\vspace{3mm}

\section{Hochschild homology}

\subsection{Graded Hecke algebras}
\label{sec:HGHA}

In \cite[Theorem 3.4]{SolHomGHA} the author showed that, for any parameter function $k$, 
\begin{equation}\label{eq:7.1}
HH_n (\mh H \rtimes \Gamma) \cong HH_n (S (\mf t^* ) \rtimes W' )
\end{equation}
as vector spaces. Here we will provide a more explicit version of \eqref{eq:7.1}
which highlights the $Z(\mh H \rtimes \Gamma)$-module structure of
$HH_* (\mh H \rtimes \Gamma)$ and is suitable for generalization to
affine Hecke algebras.

In the proof we will use some rather deep results from \cite{Sol-Irr}, which
were established under the assumption that all the parameters $k_\alpha$ are real.
In view of the scaling isomorphisms \eqref{eq:1.3} the only cases not covered
are those where $R$ contains nonperpendicular roots $\alpha, \beta$ of different
length, such that $k_\alpha$ and $k_\beta$ are $\R$-linearly independent in $\C$.

Let $\langle W' \rangle$ be a set of representatives for the conjugacy classes in
$W' = W \rtimes \Gamma$. For each $w \in \langle W' \rangle$ we pick an analytic
family of $\mh H$-representations $\{ \pi_{\rho_w} (\lambda) \mid \lambda \in t^{\rho_w} \}$
as in Corollary \ref{cor:6.5}. This yields an algebra homomorphism
\begin{equation}\label{eq:7.2}
\tilde \rho : \mh H \rtimes \Gamma \to \bigoplus_{w \in \langle W' \rangle} 
\mc O (\mf t^{\rho_w}) \otimes \mr{End}_\C (V_{\rho_w}) ,
\end{equation}
where $V_{\rho_w}$ is the vector space underlying the representations $\pi_{\rho_w} (\lambda)$.
Since these $\mh H \rtimes \Gamma$-representations are irreducible for generic 
$\lambda ,\: Z( \mh H \rtimes \Gamma) \cong \mc O (\mf t)^{W'}$ acts on them by scalars.
In other words,
\begin{equation}\label{eq:7.3}
\tilde \rho (Z (\mh H \rtimes \Gamma)) \subset \bigoplus_{w \in \langle W' \rangle} 
\mc O (\mf t^{\rho_w}) = Z \big( \bigoplus_{w \in \langle W' \rangle} 
\mc O (\mf t^{\rho_w}) \otimes \mr{End}_\C (V_{\rho_w}) \big) .
\end{equation}
Suppose that $\rho_w$ comes from a discrete series representation of $\mh H_P$. Then its
central character is an element of $\mf t_P / W (R_P) \rtimes \Gamma_\rho$, and for any
representative $cc (\rho_w) \in \mf t_P$ and any $f \in Z( \mh H \rtimes \Gamma)$: 
\begin{equation}\label{eq:7.4}
\tilde \rho (f) (w,\lambda) = f ( cc (\rho_w ) + \lambda ) .
\end{equation}
By Morita equivalence and by the Hochschild--Kostant--Rosenberg theorem
\begin{equation}\label{eq:7.5}
HH_* \big( \mc O (\mf t^{\rho_w} ) \otimes \mr{End}_\C (V_{\rho_w}) \big) \cong
HH_* \big( \mc O (\mf t^{\rho_w} ) \big) \cong \Omega^* (\mf t^{\rho_w}) ,
\end{equation}
where $\Omega^n$ stands for algebraic $n$-forms. Recall that the Morita invariance of
Hochschild homology is implemented by a generalized trace map \cite[Section I.2]{Lod}.
For $g \in \mc G_{\rho_w}$ the $\mh H \rtimes \Gamma$-representations $\pi_{\rho_w}(\lambda)$
and $\pi_{\rho_w}(g \lambda)$ have the same trace, so the image of $HH_* (\rho)$
cannot distinguish the points $(w,\lambda)$ and $(w,g \lambda)$. Thus from \eqref{eq:7.2}
and \eqref{eq:7.5} we obtain a linear map 
\begin{equation}\label{eq:7.6}
HH_* (\tilde \rho) = \bigoplus_{w \in \langle W' \rangle} HH_* (\tilde \rho )_w : 
HH_* (\mh H \rtimes \Gamma) 
\to \bigoplus_{w \in \langle W' \rangle} \Omega^* (\mf t^{\rho_w} )^{\mc G_{\rho_w}} .
\end{equation}
In view of \eqref{eq:7.4} the action of $Z (\mh H \rtimes \Gamma) \cong \mc O (\mf t)^{W'}$
on the right hand side becomes the natural evaluation if we replace $\mf t^{\rho_w}$ by 
$cc (\rho_w) + \mf t^{\rho_w}$.

\begin{thm}\label{thm:7.1}
Suppose that there exists $\ep \in \C^\times$ such that $k_\alpha \in \ep \R$ 
for all $\alpha \in R$. The map
\[
HH_* (\tilde \rho) : HH_* (\mh H \rtimes \Gamma) \to \bigoplus_{w \in \langle W' \rangle} 
\Omega^* (cc (\rho_w) + \mf t^{\rho_w} )^{\mc G_{\rho_w}}
\]
is an isomorphism of $Z (\mh H \rtimes \Gamma)$-modules.
\end{thm}
\emph{Proof.}
The isomorphism $\mh H( \tilde{\mc R},\ep k) \to \mh H (\tilde{\mc R},k)$ from \eqref{eq:1.3}
allows us to assume that $k$ is real valued. Recall that the Hochschild complex for a unital 
algebra is $C_n (A) = A^{\otimes n+1}$ with differential
\begin{multline*}
b_n (a_0 \otimes \cdots \otimes a_n) = \sum\nolimits_{i=0}^{n-1} (-1)^n a_0 
\otimes \cdots \otimes a_i a_{i+1} \otimes \cdots \otimes a_n + \\ 
(-1)^n a_n a_0 \otimes a_1 \otimes \cdots \otimes a_{n-1} .
\end{multline*}
As discussed in the proof of \cite[Theorem 3.4]{SolHomGHA}, we can choose 
representatives $w$ for the conjugacy classes in $W'$, such that the following hold:
\begin{itemize}
\item $w = \gamma w_P$ with $\gamma \in \Gamma$ and $w_P \in W (R_P)$;
\item $\mf a^w \subset \mf a^P \cap \mf a^\gamma \cap \mf a^{w_P}$;
\item all elements of $\mc O (\mf t_w) \subset \mc O (\mf t / \mf t_P) 
\subset \mh H$ commute with $w \in \mh H \rtimes \Gamma$, where $\mf t_w := 
\mf t / (w-1)\mf t \cong \mf t^w$;
\item every element of $HH_* (\mh H \rtimes \Gamma)$ can be represented by a Hochschild 
cycle in $\bigoplus_{w \in \langle W' \rangle} w C_* (\mc O (\mf t_w))$;
\item the differential complex $w C_* (\mc O (\mf t_w))$ contributes precisely
$\Omega^* (\mf t^w)^{Z_{W'}(w)}$ to\\
$HH_* (\mh H \rtimes \Gamma)$, via the natural map
\end{itemize}
\begin{equation}\label{eq:7.11}
C_* (\mc O (\mf t_w)) \to \Omega^* (\mf t^w) : a_0 \otimes a_1 \otimes \cdots \otimes a_n 
\mapsto a_0 d a_1 \cdots d a_n .
\end{equation}
In the last point it also suffices to take a subspace (say $w \Omega^*_w$) of the cycles in 
$w C_* (\mc O (\mf t_w))$ which is mapped bijectively to $\Omega^* (\mf t^w)^{Z_{W'}(w)}$
by \eqref{eq:7.11}. We intend to show that the map 
\begin{equation}\label{eq:7.13}
\bigoplus_{w \in \langle W' \rangle} w \Omega_w^* \to \bigoplus_{w \in \langle W' \rangle} 
\Omega^* (\mf t^{\rho_w} )^{\mc G_{\rho_w}}
\end{equation}
induced by $HH_* (\tilde \rho)$ is bijective. Our strategy is first to show that this map is 
upper triangular with respect to some reflexive transitive order on $\langle W' \rangle$ 
and then to apply the results of Section \ref{sec:complex} to the 
elements of $\langle W' \rangle$ that equivalent in this sense.

For the first point we study what happens to $w C_* (\mc O (\mf t_w))$ under 
$HH_* (\tilde \rho )_{w'}$. We need to apply the generalized trace map
\[
HH_* (\mr{gtr}) : HH_* \big( \mc O (\mf t^{\rho_{w'}}) \otimes \mr{End}_\C (V_{\rho_{w'}}) \big) 
\to HH_* (\mc O (\mf t^{\rho_{w'}})) ,
\]
which is defined on the level of Hochschild complexes by
\[
\mr{gtr} (m_0 a_0 \otimes m_1 a_1 \otimes \cdots \otimes m_n a_n ) = 
\mr{tr}(m_0 m_1 \cdots m_n) a_0 \otimes a_1 \otimes \cdots \otimes a_n
\]
for endomorphisms $m_i \in \mr{End}_\C (V_{\rho_w})$ and functions 
$a_i \in \mc O (\mf t^{\rho_{w'}})$. Since $w C_* (\mc O (\mf t_w))$ is contained in the 
Hochschild complex of the commutative algebra 
\[
A_w = \C [C_w ] \otimes \mc O (\mf t_w)
\]
generated by $w$ and $\mc O (\mf t_w)$, we can analyse 
$HH_* (\tilde \rho )_{w'} \big(w C_* (\mc O (\mf t_w)) \big)$ via a filtration of the
representation $\pi_{\rho_{w'}} (\lambda)$ by onedimensional $A_w$-modules. With respect to
a corresponding basis of $V_{\rho_{w'}}$, all the matrices $\pi_{\rho_{w'}} (\lambda) (A_w)$ 
are upper triangular. Combining that with the formula for the generalized trace map, one 
concludes that $HH_* (\tilde \rho )_{w'} (A_w)$ sees $\pi_{\rho_{w'}} (\lambda)$ only 
via the diagonal parts of these matrices, that is, via the semisimplication of this 
$A_w$-representation.
 
Let us assume for the moment that $\lambda \in \mf t_w$ is generic.
The semisimplification of the restriction of $\pi_{\rho_{w'}} (\lambda)$ to 
$\mc O (\mf t_w)$ is a direct sum of onedimensional modules with characters $s + g(\lambda)$, 
where $s$ runs through the $\mc O (\mf t_w)$-characters of $\rho_{w'}$ and $g$ runs through 
$\mc W'_{\rho_{w'}}$. As $\lambda$ is generic, this yields a decomposition into a direct 
sum of exactly $|\mc W'_{\rho_{w'}}| \; \mc O (\mf t_w)$-submodules, the summand 
parametrized by $g \in \mc W'_{\rho_{w'}}$ having weights $s + g(\lambda)$.
The intertwining operator associated to $g$, as in \cite[Section 8]{SolGHA}, provides a linear
bijection between the summands corresponding to $(\lambda,1)$ and to $(g^{-1} \lambda, g)$. 

From the above submodules of the semisimplification of $\pi_{\rho_{w'}} (\lambda)$ we take 
the ones with the same $\mc O (\mf t_w)$-character together, thus creating a decomposition into 
$\mc O (\mf t_w)$-submodules $M_i (\lambda)$ with pairwise different weights. The $M_i (\lambda)$ 
are stabilized by $w$, because $w$ commutes with $\mc O (\mf t_w)$. As a representation 
of the finite group $C_w$ generated by $w ,\; M_i (\lambda)$ is independent of
$\lambda$ because $\lambda$ can be varied continuously. Neither does it depend on $i$, since
$M_i (\lambda)$ is isomorphic to $M_j (g \lambda)$ via the intertwining operator associated to 
a suitable $g \in \mc W'_{\rho_{w'}}$. We conclude that the semisimplification of 
$\pi_{\rho_{w'}} (\lambda)$ as an $A_w$-module is of the form
\begin{equation}\label{eq:7.7}
V_w \otimes V_{\mc O} (\lambda)
\end{equation}
for some $C_w$-representation $V_w$ and some $\mc O (\mf t_w)$-module $V_{\mc O} (\lambda)$. 
Although the above argument uses that $\lambda$ is generic, the conclusion extends to all 
$\lambda \in \mf t_w$ because we are dealing with a continuous family of representations.

From \eqref{eq:7.7} we see that for any $f_i \in \mf t_w$
\begin{equation}\label{eq:7.8}
HH_n (\tilde \rho )_{w'} (w f_0 \otimes f_1 \cdots \otimes f_n) = 
\mr{tr}(w,V_w) \dim (V_w )^{-1} HH_* (\tilde \rho )_{w'} (f_0 \otimes f_1 \cdots \otimes f_n) .
\end{equation}
If $(w-1) \mf t$ is not contained in $g (w'-1) \mf t$ for any $g \in W'$ and 
$\lambda \in \mf t_{w'}$ is generic, then $w$ permutes the $\mc O (t_w)$-weight spaces of 
$\pi_{\rho_{w'}}(\lambda)$ in a nontrivial way. Hence tr$(w,V_w) = 0$ in such cases, and by 
continuity this extends to all $\lambda \in \mf t_w$. For reference we state this as
\begin{equation}\label{eq:7.9}
HH_* (\tilde \rho )_{w'} \big( w C_* (\mc O (t_w)) \big) = 0 
\quad \text{when} \quad (w-1) \mf t \not\subset W'(w'-1) \mf t_w .
\end{equation}
Now pick $\lambda \in i \mf a^{w'}$ with $\norm{\lambda}$ so small that 
$g(\lambda) = \lambda$ is equivalent to $\mf t^g \supset \mf t^{w'}$ for $g \in W'$.
By parts (3) and (4) of Theorem \ref{thm:6.4}, among the $\mc O (\mf t) \rtimes 
W'$-representations $\zeta^\vee (\pi_{\rho_g} (\mu))$, only those with $\mu \in W' \lambda$ 
and $W' \mf t^{\rho_g} \supset \mf t^{w'}$ have central character $W' \lambda$. For the 
duration of the proof we write $g \geq w'$ to describe this situation. Similarly we write 
$g > w'$ if $g \geq w'$ but $W' \mf t^{\rho_g} \neq W' \mf t^{w'}$, and we write 
$g \approx w'$ if $W' \mf t^{\rho_g} = W' \mf t^{w'}$. To simplify the situation somewhat, 
we may and will assume that the elements of $\langle W' \rangle$ are chosen such that 
$\mf t^{\rho_g} \supset \mf t^{w'}$ for all $g \in \langle W' \rangle_{\geq w'}$.

By \eqref{eq:6.9} the number of irreducible $\mc O (\mf t) \rtimes W'$-representations with 
central character $W' \lambda$ is the number of conjugacy classes in the isotropy group 
$W'_\lambda$. Together with the bijectivity of $\zeta^\vee$ this shows that every conjugacy 
class of $W'_\lambda$ is contained in a unique conjugacy class of $W'$, and that all elements 
of $\langle W' \rangle_{\geq w'}$ are obtained in this way. Given $z \in Z_{W'}(w')$,
we have $z \lambda \in \mf t^{w'} \subset \mf t^{\rho_g}$ for all $g \in \langle W' 
\rangle_{\geq w'}$. By Lemma \ref{lem:6.3}.c and Corollary \ref{cor:6.5} the representations 
\[
\{ \zeta^\vee (\pi_{\rho_g}(z \lambda)) : g \in \langle W' \rangle_{\geq w'} \}
\]
form a basis of the subspace of $G_\Q (\mh H \rtimes \Gamma)$ corresponding to the central 
character $W' \lambda$. By parts (a) and (b) of Lemma \ref{lem:6.3} this means that 
$\pi_{\rho_g}(z \lambda)$ and $\pi_{\rho_g}(\lambda)$ are $\mc W'_{\rho_g}$-associate,
for all $g \in \langle W' \rangle_{\geq w'}$. For similar reasons other generic $\lambda'
\in \mf t^{w'}$ yield representations with other traces. For $g \approx w'$ this allows us 
to conclude that the orbifolds $\mf t^{\rho_g} / \mc W'_{\rho_g}$ and 
$\mf t^{w'} / Z_{W'}(w')$ have the same generic parts, and hence are isomorphic.

Recall from Theorem \ref{thm:6.4} that $\zeta^\vee$ does not change the $W'$-type of 
representations and induces a bijection between Grothendieck groups. Therefore the matrix 
\[
\big( \mr{tr}(w, \pi_{\rho_g} (\lambda)) \big)_{w,g \in \langle W' \rangle_{\geq w'}}
\]
is invertible. We note that this matrix does not depend on $\lambda \in \mf t^w$.
As discussed before \eqref{eq:7.9}, all entries with $w > g$ are 0, so the submatrix 
$\big( \mr{tr}(w, \pi_{\rho_g} (\lambda)) \big)_{w,g \in \langle W' \rangle_{\approx w'}}$
is also invertible. Together with Corollary \ref{cor:6.5} and \eqref{eq:7.8} this implies 
that the map
\begin{equation}\label{eq:7.10}
\bigoplus_{w \in \langle W' \rangle_{\approx w'}} w \mc O (\mf t^w)^{Z_{W'}(w)} \to
\bigoplus_{g \in \langle W' \rangle_{\approx w'}} \mc O (\mf t^{\rho_g})^{\mc W'_{\rho_g}}
\end{equation}
induced by $HH_0 (\tilde \rho)$ becomes a bijection upon specializing at any 
$\lambda \in \mf t^{w'}$. As shown above, the underlying orbifolds on both sides are 
$\langle W' \rangle_{\approx w'} \times (\mf t^w)^{Z_{W'}(w)}$, so by \eqref{eq:7.8} 
we can regard \eqref{eq:7.10} as an invertible complex matrix of size 
$| \langle W' \rangle_{\approx w'} |$ tensored with the unique isomorphisms of orbifolds 
$\mf t^{\rho_g} / \mc W'_{\rho_g} \to \mf t^w / Z_{W'}(w')$ that respect the projections 
to $\mf t / W'$. This description and \eqref{eq:7.8} show that the natural extension of 
\eqref{eq:7.10} to differential forms,
\begin{equation}\label{eq:7.12}
\bigoplus_{w \in \langle W' \rangle_{\approx w'}} w \Omega^* (\mf t^w)^{Z_{W'}(w)} \to
\bigoplus_{g \in \langle W' \rangle_{\approx w'}} \Omega^* (\mf t^{\rho_g})^{\mc W'_{\rho_g}} ,
\end{equation}
is also bijective. In view of \eqref{eq:7.11} the map \eqref{eq:7.12} is induced by 
$HH_* (\tilde \rho)$. Together with \eqref{eq:7.9} this proves that \eqref{eq:7.13} 
is bijective, as required. $\qquad \Box$ 
\vspace{2mm}

In fact the above proof establishes a slightly more precise statement: 
the map of differential complexes
\begin{equation}\label{eq:7.16}
\text{gtr} \circ C_* (\tilde \sigma) : C_* (\mh H \rtimes \Gamma) \to
\bigoplus_{w \in \langle W' \rangle} C_* (\mf t^{\rho_w} )^{{\mc W'}_{\rho_w}}
\end{equation}
is a quasi-isomorphism.

Theorem \ref{thm:7.1} has a direct analogue for the (periodic) cyclic homology of graded
Hecke algebras:

\begin{cor}\label{cor:7.2}
Let $k$ be as in Theorem \ref{thm:7.1}.
The algebra homomorphism  \eqref{eq:7.2} induces isomorphisms
\[
\begin{array}{llll}
HC_* (\tilde \rho) : & HC_* (\mh H \rtimes \Gamma) & \to &  
\bigoplus\limits_{w \in \langle W' \rangle} HC_* ( \mc O (\mf t^{\rho_w}) )^{\mc W_{\rho_w}} , \\
HP_{\text{ev/odd}} (\tilde \rho) : & HP_{\text{ev/odd}} (\mh H \rtimes \Gamma) & \to &  
\bigoplus\limits_{w \in \langle W' \rangle} H_{DR}^{\text{ev/odd}} (\mf t^{\rho_w} 
)^{\mc W_{\rho_w}} \; \cong \; H_{DR}^{\text{ev/odd}} (\tilde{\mf t})^{W'} .
\end{array} 
\]
\end{cor}
We remark that the (periodic) cyclic homology of smooth commutative algebras like
$\mc O (\mf t^{\rho_w})$ is well-known from Hochschild--Kostant--Rosenberg theorem, 
see for example \cite[Theorem 3.4.12]{Lod}. 

\emph{Proof.}
Recall that the relation between the Hochschild complex $C_* (A)$ and the cyclic bicomplex 
$CC_{**}(A)$ of an algebra $A$ gives rise to Connes' periodicity exact sequence:
\begin{equation}\label{eq:7.15}
\cdots \to HH_n (A) \to HC_n (A) \to HC_{n-2}(A) \to HH_{n-1}(A) \to \cdots  
\end{equation}
In particular this holds for $A = \mh H \rtimes \Gamma$. The algebra homomorphism $\tilde \rho$ 
and the generalized trace map send \eqref{eq:7.15} to the corresponding long exact sequence 
for the cyclic bicomplex 
$\bigoplus_{w \in \langle W' \rangle} CC_{**} (\mc O (\mf t^{\rho_w}) )^{\mc W_{\rho_w}}$.
Theorem \ref{thm:7.1} says that this map induces an isomorphism on Hochschild homology. Since 
$HC_0 = HH_0$, it follows with induction to $n$ from \eqref{eq:7.15} and the five lemma that
$\tilde \rho$ induces the asserted isomorphism on $HC_*$. Since $HH_n (\mh H \rtimes \Gamma) = 0$
for $n > \dim_\C \mf t$, it follows that $HC_{n+2} (\mh H \rtimes \Gamma) \cong HC_n  
(\mh H \rtimes \Gamma)$ for such $n$. By \cite[Proposition 5.1.9]{Lod} this implies that 
$HP_n (\mh H \rtimes \Gamma) \cong HC_n  (\mh H \rtimes \Gamma)$ for all $n > \dim_\C \mf t$. 
The same holds for the (periodic) cyclic homology of $\mc O (\mf t^{\rho_w})$, 
so $\tilde \rho$ also induces an isomorphism 
\begin{equation}\label{eq:7.14}
HP_{\text{ev/odd}} (\mh H \rtimes \Gamma) \to \bigoplus\limits_{w \in \langle W' \rangle} 
HP_{\text{ev/odd}} (\mc O (\mf t^{\rho_w} ) )^{\mc W_{\rho_w}} \cong 
\bigoplus\limits_{w \in \langle W' \rangle} 
H_{DR}^{\text{ev/odd}} (\mf t^{\rho_w} )^{\mc W_{\rho_w}} .
\end{equation}
By Corollary \ref{cor:6.5} and \eqref{eq:6.10}
\[
\bigsqcup_{w \in \langle W' \rangle} \mf t^{\rho_w} / \mc W_{\rho_w} \cong  
\bigsqcup_{w \in \langle W' \rangle} \mf t^w / Z_{W'}(w) \cong \tilde{\mf t} / W'
\]
as orbifolds. Therefore the right hand side of \eqref{eq:7.14} is isomorphic to
$H_{DR}^{\text{ev/odd}} (\tilde{\mf t})^{W'}. \qquad \Box$
\vspace{2mm}

\subsection{Affine Hecke algebras}
\label{sec:Hhaha} 

We will transfer the desciption of the Hochschild homology of graded Hecke algebras,
as given in the previous section, to affine Hecke algebras. The connection is provided by
Lusztig's reduction theorems. Since these involve formal completions of algebras, we spend
a few words on the modifications which are necessary for homological algebra in that context.

Given an algebra $A$ and a central ideal $I$ of the multiplier algebra of $A$, we endow
\[
\widehat{A} := \lim_{\infty \leftarrow n} A / I^n A 
\]
with the $I$-adic topology. Then it is natural to use completed tensor products and relative
homological algebra, as in \cite[Section 2]{KNS}. The resulting Hochschild homology is
denoted $HH_*^{top} (\widehat A)$ and it is a module over $\widehat{Z(A)}$.

Let $V$ be any complex affine variety. Recall that a finite type $\mc O (V)$-algebra is an
algebra $A$ together with a homomorphism from $\mc O (V)$ to the centre of the
multiplier algebra of $A$, which makes $A$ into a $\mc O (V)$-module of finite rank.

\begin{thm}\label{thm:8.1} 
\textup{\cite[Theorem 3]{KNS}} \\
Let $A$ be a finite type $\mc O (V)$-algebra and let $I$ be an ideal of $\mc O (V)$. 
The inclusion $A \to \widehat{A}$ induces an isomorphism
\[
HH_* (A) \otimes_{\mc O (V)} \widehat{\mc O (V)} \to HH_*^{top} (\widehat{A}) . 
\]
\end{thm}

Because we intend to use the results from Section \ref{sec:complex},
we fix $\ep \in \C \setminus i \R$ and we assume that $q_{\alpha^\vee} \in \exp (\ep \R)$
for every $\alpha^\vee \in \R_{nr}^\vee$ 
We pick analytic families of $\mc H \rtimes \Gamma$-representations as in \eqref{eq:6.3} and
Theorem \ref{thm:6.6}. These give rise to an algebra homomorphism 
\begin{equation}\label{eq:8.1}
\tilde \sigma : \mc H \rtimes \Gamma \to 
\bigoplus_{w \in \langle W' \rangle} \bigoplus_{i=1}^{c(w)}
\mc O (T^{\sigma_{w,i}}) \otimes \mr{End}_\C (V_{\sigma_{w,i}}) . 
\end{equation}
This induces a map $C_* (\tilde \sigma)$ between the corresponding Hochschild complexes, 
which we can compose with the generalized trace map to reach the direct sum of the 
Hochschild complexes $C_* (\mc O (T^{\sigma_{w,i}}))$. Since the image cannot distinguish 
representations with the same trace, we obtained a morphism of differential complexes
\begin{equation}\label{eq:8.3}
\text{gtr} \circ C_* (\tilde \sigma) : C_* (\mc H \rtimes \Gamma) \to
\bigoplus_{w \in \langle W' \rangle} \bigoplus_{i=1}^{c(w)}
C_* \big( \mc O (T^{\sigma_{w,i}}) \big)^{\mc G_{\sigma_{w,i}}}
\end{equation}
Just as for graded Hecke algebras the induced map on Hochschild homology lands in a space
of invariant differential forms:
\begin{equation}\label{eq:8.2}
HH_* (\tilde \sigma) : HH_* (\mc H \rtimes \Gamma) \to \bigoplus_{w \in \langle W' \rangle} 
\bigoplus_{i=1}^{c(w)} \Omega^* (cc (\sigma_{w,i}) T^{\sigma_{w,i}} )^{\mc G_{\sigma_{w,i}}} .
\end{equation}
Here $cc (\sigma_{w,i})$ is an $\mc A$-weight of the irreducible $\mc H^P \rtimes 
\Gamma_{\sigma_{w,i}}$-representation $\sigma_{w,i}$, which serves mainly to make the 
module structure over the centre clear. An element 
$f \in Z (\mc H \rtimes \Gamma) \cong \mc O (T )^{W'}$ acts on the right hand side of 
\eqref{eq:8.2} via evaluation on $cc (\sigma_{w,i}) T^{\sigma_{w,i}} \subset T$. 
This construction makes \eqref{eq:8.2} $Z(\mc H \rtimes \Gamma)$-linear.

\begin{thm}\label{thm:8.2}
Suppose that there exists $\ep \in \C \setminus i \R$ such that $q_{\alpha^\vee} \in \exp (\ep \R)$
for all $\alpha^\vee \in R_{nr}^\vee$. The map \eqref{eq:8.3} is a quasi-isomorphism and 
\eqref{eq:8.2} is an isomorphism of $Z( \mc H \rtimes \Gamma)$-modules. 
\end{thm}
\emph{Proof.}
We only have to deal with the second claim, since it encompasses the first.
Pick any $u \in T_{un}$ and $\lambda \in \ep \mf a$. By Theorem \ref{thm:2.2} the algebras
$\mc H \rtimes \Gamma$ and $\mh H (\tilde{\mc R}^{P(u)},k^{P(u)}) \rtimes W'_{P(u)}$ become 
isomorphic upon formally completing at $W' u \exp (\lambda)$, respectively at 
$W'_{\Z P(u)} \lambda$. This enables us to transform the $\pi_{\sigma_{w,i}}(t)$ for which 
\[
cc (\sigma_{w,i}) T^{\sigma_{w,i}} \cap W' u \exp (\ep \mf a) \neq \emptyset
\]
into analytic families of $\mh H (\tilde{\mc R}^{P(u)},k^{P(u)}) \rtimes W'_{P(u)}$-representations
\[
\{ \pi_{\rho_w}(\lambda) : w \in \langle W'_u \rangle , \lambda \in \mf t^{\rho_w} \} ,
\]
like we did in the proof of Theorem \ref{thm:6.6}. By Theorem \ref{thm:6.6} the formal completions
\begin{align*}
& \bigoplus_{w \in \langle W' \rangle} \bigoplus_{i=1}^{c(w)} 
\widehat{\mc O ( cc (\sigma_{w,i}) T^{\sigma_{w,i}})^{\mc G_{\sigma_{w,i}}} }_{W' u \exp (\lambda)} 
\otimes \mr{End}_\C (V_{\sigma_{w,i}}) \\
& \text{and} \quad \bigoplus_{w \in \langle W'_u \rangle} 
\widehat{\mc O (cc (\rho_w) + \mf t^{\rho_w})^{\mc W_{\rho_w}}}_{W'_{\Z P(u)} \lambda} 
\otimes \mr{End}_\C (V_{\rho_w})
\end{align*}
are Morita equivalent. All this makes it possible to fit \eqref{eq:8.2} in a commutative diagram
\[
\begin{array}{ccc}
HH_* (\mc H \rtimes \Gamma) & \to & \bigoplus\limits_{w \in \langle W' \rangle} 
\bigoplus\limits_{i=1}^{c(w)} \Omega^* (cc (\sigma_{w,i}) T^{\sigma_{w,i}} )^{\mc G_{\sigma_{w,i}}} \\
\downarrow & & \downarrow \\
HH_* \big( \widehat{\mc H}_{W' u \exp (\lambda)} \rtimes \Gamma \big) & \to &
\hspace{-3mm} \bigoplus\limits_{w \in \langle W' \rangle} \bigoplus\limits_{i=1}^{c(w)} 
\widehat{\Omega^* (cc (\sigma_{w,i}) T^{\sigma_{w,i}} )^{\mc G_{\sigma_{w,i}}}}_{W' u \exp (\lambda)} \\
\downarrow & & \downarrow \\
\hspace{-2mm} HH_* \big( \widehat{\mh H (\tilde{\mc R}^{P(u)}, k^{P(u)})_{W'_{\Z P(u)} \lambda}} 
\rtimes W'_{P(u)} \big) \hspace{-3mm} & \to & \bigoplus\limits_{w \in \langle W'_u \rangle} 
\widehat{\Omega^* (cc (\rho_w) + \mf t^{\rho_w})^{\mc W_{\rho_w}} }_{W'_{\Z P(u)} \lambda}
\end{array} 
\]
By Theorems \ref{thm:7.1} and \ref{thm:8.1} the lower horizontal map is an isomorphism of \\ 
$Z \big( \widehat{\mh H (\tilde{\mc R}^{P(u)}, 
k^{P(u)})_{W'_{\Z P(u)} \lambda}} \rtimes W'_{P(u)} \big)$-modules.
By the above constructions the same goes for the two lower vertical maps, and hence for the middle 
horizontal map as well. Thus the diagram and Theorem \ref{thm:8.1} show that \eqref{eq:8.2} 
is an homomorphism of $Z (\mc H \rtimes \Gamma)$-modules which becomes an isomorphism upon 
completing with respect to an arbitrary maximal ideal of this commutative algebra. 

Now consider any $x \in \ker HH_* (\tilde \sigma)$. Then 
\[
Z (\mc H \rtimes \Gamma) x \cong Z (\mc H \rtimes \Gamma) / J_x
\]
for some ideal $J_x \subset Z (\mc H \rtimes \Gamma)$. The above implies that 
$\widehat{Z (\mc H \rtimes \Gamma) / J_x} = 0$ for every maximal ideal of 
$Z (\mc H \rtimes \Gamma)$, which means that $J_x$ is not contained in any maximal ideal. 
Hence $J_x = Z (\mc H \rtimes \Gamma)$ and $x = 0$. A similar argument shows that coker 
$HH_* (\tilde \sigma) = 0$, so \eqref{eq:8.2} is an isomorphism. $\qquad \Box$ 
\vspace{2mm}

\begin{cor}\label{cor:8.3}
The algebra homomorphism  \eqref{eq:8.2} induces isomorphisms
\[
\begin{array}{llll}
\hspace{-2mm} HC_* (\tilde \sigma) : & HC_* (\mc H \rtimes \Gamma) & \!\!\to\!\! &  
\bigoplus\limits_{w \in \langle W' \rangle} \bigoplus\limits_{i=1}^{c(w)} 
HC_* ( \mc O (T^{\sigma_{w,i}} ))^{\mc G_{\sigma_{w,i}}} , \\
\hspace{-2mm} HP_{\text{ev/odd}} (\tilde \sigma) : \vspace{-2mm} & 
HP_{\text{ev/odd}} (\mc H \rtimes \Gamma) & \!\!\to\!\! &  
\bigoplus\limits_{w \in \langle W' \rangle} \bigoplus\limits_{i=1}^{c(w)} H_{DR}^{\text{ev/odd}} 
(T^{\sigma_{w,i}} )^{\mc G_{\sigma_{w,i}}}  \cong  H_{DR}^{\text{ev/odd}} (\tilde T)^{W'} .
\end{array} 
\]
\end{cor}
\emph{Proof.}
This follows in the same way from Theorem \ref{thm:8.2} as Corollary \ref{cor:7.2} follows from 
Theorem \ref{thm:7.1}. $\qquad \Box$ 
\vspace{2mm}

We note that the periodic cyclic homology of $\mc H (\mc R,q) \rtimes \Gamma$ was already 
known from \cite[Section 5.2]{Sol-Irr} for positive $q$. However, there it was obtained in a 
rather indirect way, via the topological $K$-theory of a $C^*$-completion of 
$\mc H (\mc R,q) \rtimes \Gamma$. With Corollary \ref{cor:8.3} we can understand 
$HP_* (\mc H (\mc R,q) \rtimes \Gamma)$ more explicitly and for more general parameter 
functions $q$.

Maybe one can improve on Theorem \ref{thm:8.2} by finding suitable maps 
$\phi_{w,i} : \mc H \rtimes \Gamma \to A_{w,i}$, where $A_{w,i}$ is an algebra with dual 
space $T^w_i / Z_{W'}(w)$ and $HH_* (A_{w,i}) \cong \Omega^* (T^w_i )^{Z_{W'}(w)}$. 
Yet this is difficult in the present setup, since the representations $\pi_{\sigma_{w,i}}(t)$ 
and $\pi_{\sigma_{w,i}}(g t)$ with $g \in \mc G_{\sigma_{w,i}}$ need not be isomorphic. 

For equal parameter functions $q$, Lusztig achieved something similar with an "asymptotic Hecke
algebra" \cite{Lus-C2}. This was used in \cite{BaNi} to compute the periodic cyclic homology of
$\mc H (\mc R,q)$. It is not known whether this strategy can be employed for Hochschild homology
or for unequal parameter functions.
\vspace{2mm}

\subsection{Schwartz algebras}
\label{sec:homS}

In this section we will compute the Hochschild homology of the Schwartz completion 
$\mc S \rtimes \Gamma$ of an (extended) affine Hecke algebra $\mc H \rtimes \Gamma$ with 
positive parameters $q$. It turns out that there is a clear relation between 
$HH_* (\mc S \rtimes \Gamma)$ and $HH_* (\mc H \rtimes \Gamma)$, similar to the relation 
between $C^\infty (S^1)$ and $\mc O (\C^\times ) \cong \C [z,z^{-1}]$. 

To get the correct answer we must take the topology of $\mc S \rtimes \Gamma$ into account. 
The best way to do so is with complete bornological algebras \cite{Mey-Emb} and bornological 
tensor products, which we will denote by $\widehat \otimes$. Thus 
$HH_* (\mc S \rtimes \Gamma)$ is the homology of the differential complex 
$\big( (\mc S \rtimes \Gamma)^{\widehat \otimes n + 1}, b_n \big)$. It is isomorphic to 
\[
\mr{Tor}_*^{\mc S \rtimes \Gamma \widehat \otimes (\mc S \rtimes \Gamma)^{op}}
(\mc S \rtimes \Gamma, \mc S \rtimes \Gamma) ,
\]
computed in the category of complete bornological $\mc S \rtimes \Gamma$-bimodules.

We will use the notation from \eqref{eq:6.3}. By Theorem \ref{thm:5.2} the evaluation map 
\[
\mc H \rtimes \Gamma \to \mc O (T^{\sigma_{w,i}}) \otimes \mr{End}_\C (V_{\sigma_{w,i}})
\]
extends to a homomorphism of Fr\'echet algebras
\[
\mc S \rtimes \Gamma \to 
C^\infty (T_{un}^{\sigma_{w,i}}) \otimes \mr{End}_\C (V_{\sigma_{w,i}}) .
\]
Applying this to all $w$ and $i$ yields a homomorphism
\begin{equation}\label{eq:9.6}
\tilde \sigma_{un} : \mc S \rtimes \Gamma \to \bigoplus_{w \in \langle W' \rangle} 
\bigoplus_{i=1}^{c(w)} C^\infty (T_{un}^{\sigma_{w,i}}) \otimes \mr{End}_\C (V_{\sigma_{w,i}})
\end{equation}
which extends \eqref{eq:8.1}. Recall that by the topological version of the 
Hochschild--Kostant--Rosenberg Theorem 
\[
HH_* (C^\infty (T_{un})) \cong \Omega_{sm}^* (T_{un}) ,
\]
the algebra of smooth complex valued differential forms on the real manifold $T_{un}$.
Hence \eqref{eq:9.6}, like $\tilde \sigma$, induces a $\C$-linear map
\begin{equation} \label{eq:9.5} 
HH_* (\tilde \sigma_{un}) : HH_* (\mc S \rtimes \Gamma) \to \bigoplus_{w \in \langle W' \rangle} 
\bigoplus_{i=1}^{c(w)} \Omega_{sm}^* (T_{un}^{\sigma_{w,i}})^{\mc G_{\sigma_{w,i}}} .
\end{equation}
It becomes $Z (\mc S \rtimes \Gamma)$-linear if we write
$cc (\sigma_{w,i}) T^{\sigma_{w,i}}_{un}$ instead of $T_{un}^{\sigma_{w,i}}$ and 
interpret the $Z (\mc S \rtimes \Gamma)$-action accordingly.

\begin{thm} \label{thm:9.1} 
Let $q$ be positive. The map \eqref{eq:9.5} is an isomorphism of 
$Z( \mc S \rtimes \Gamma)$-modules and it extends Theorem \ref{thm:8.2}.
\end{thm}
\emph{Proof.}
Since $e_{P,\delta} \in \mc S \rtimes \Gamma$ is a central idempotent,
\begin{equation} \label{eq:9.7}
HH_n (e_{P,\delta} \mc S \rtimes \Gamma) = 
H_n \big( (e_{P,\delta} \mc S \rtimes \Gamma)^{\widehat \otimes (*+1)}, b_* \big) \cong
H_n \big( e_{P,\delta} \mc S \rtimes \Gamma \widehat \otimes 
\mc S \rtimes \Gamma)^{\widehat \otimes *}, b_* \big) .
\end{equation}
For $\xi = (P,\delta,t) \in \Xi_{un}$ let $\mc F_\xi$ be the ring of formal power series on
$(P,\delta,T^P_{un})$ centred at $\xi$. By \cite[Theorem 3.5]{OpSo3} the functor
\[
\mc F_\xi^{\mc G_{P,\delta}} \widehat \otimes_{Z (\mc S \rtimes \Gamma)} = 
\mc F_\xi^{\mc G_{P,\delta}} \widehat \otimes_{e_{P,\delta} Z (\mc S \rtimes \Gamma)}
\]
is exact on a large class of $Z(\mc S \rtimes \Gamma)$-modules, which contains all modules 
that we use here. Applying this to the right-hand side of \eqref{eq:9.7} we obtain
\begin{equation} \label{eq:9.8}
\mc F_\xi^{\mc G_{P,\delta}} \widehat \otimes_{Z (\mc S \rtimes \Gamma)}
HH_n (e_{P,\delta} \mc S \rtimes \Gamma ) \cong 
H_n \big( \mc F_\xi^{\mc G_{P,\delta}} \widehat \otimes_{Z (\mc S \rtimes \Gamma)}
e_{P,\delta} \mc S \rtimes \Gamma \widehat \otimes 
\mc S \rtimes \Gamma)^{\widehat \otimes *}, b_* \big) .
\end{equation}
It follows from a result of Borel \cite[Remarque IV.3.5]{Tou} that the Taylor series map
\[
e_{P,\delta} Z (\mc S \rtimes \Gamma) \cong C^\infty (T^P_{un})^{\mc G_{P,\delta}} \to
\mc F_\xi^{\mc G_{P,\delta}}
\]
is surjective, and (by definition) its kernel $I_\xi^\infty$ consists of all functions
that are flat at $\xi$. Clearly the $Z(\mc S \rtimes \Gamma)$-action on \eqref{eq:9.8} 
factors through $\mc F_\xi^{\mc G_{P,\delta}}$, so $I_\xi^\infty$ annihilates \eqref{eq:9.8}. 
Together with Lemma \ref{lem:5.3} this shows that \eqref{eq:9.8} is isomorphic to 
\begin{multline}\label{eq:9.10}
\hspace{-3mm} H_n \big( (\mc F_\xi^{\mc G_{P,\delta}} \widehat \otimes_{Z (\mc S \rtimes \Gamma )}
e_{P,\delta} \mc S \rtimes \Gamma )^{\widehat \otimes (*+1)}, b_* \big) \cong
H_n \big( (\mc F_\xi^{\mc G_{P,\delta}} \otimes_{e_{P,\delta} Z (\mc H \rtimes \Gamma )}
e_{P,\delta} \mc H \rtimes \Gamma )^{\widehat \otimes (*+1)}, b_* \big) \\
= HH_n \big( \mc F_\xi^{\mc G_{P,\delta}} \otimes_{e_{P,\delta} Z (\mc H \rtimes \Gamma )}
e_{P,\delta} \mc H \rtimes \Gamma \big) . 
\end{multline}
Let $pr : \Xi \to T / W'$ be the projection that sends $\xi$ to the $Z (\mc H \rtimes 
\Gamma)$-character of the representation $\pi^\Gamma (\xi)$. 
With \eqref{eq:5.2} and Lemma \ref{lem:6.3}.c we can associate to every tempered analytic 
family of $\mc H \rtimes \Gamma$-representations a unique $(P,\delta) \in \mc P$. 
In the notation of \eqref{eq:6.3} we write $\sigma_{w,i} \prec (P,\delta)$ if 
$\pi_{\sigma_{w,i}}(t) (e_{P,\delta}) \neq 0$ for any $t \in T^P_{un}$. 

If we want to apply the algebra 
homomorphism $\tilde \sigma$ (or $\tilde \sigma_{un}$, that comes down to the same thing) to 
$\mc F_\xi^{\mc G_{P,\delta}} \otimes_{e_{P,\delta} Z (\mc H \rtimes \Gamma )} e_{P,\delta} 
\mc H \rtimes \Gamma$, we must replace the target by
\begin{multline*}
\mc F_\xi^{\mc G_{P,\delta}} \otimes_{e_{P,\delta} Z (\mc H \rtimes \Gamma )}
\bigoplus_{w \in \langle W' \rangle} \bigoplus_{i=1}^{c(w)} \mc O (T^{\sigma_{w,i}}) 
\otimes \mr{End}_\C (V_{\sigma_{w,i}}) \\
= \bigoplus_{\sigma_{w,i} \prec (P,\delta)} \widehat{\mc O (T^{\sigma_{w,i}})}_{pr (\xi)} 
\otimes \mr{End}_\C (V_{\sigma_{w,i}}) .
\end{multline*}
Let us denote the resulting algebra homomorphism by $\tilde \sigma_\xi$. It induces a morphism
$C_* (\tilde \sigma_\xi)$ on Hochschild complexes, whose composition with the generalized
trace map gtr lands in a sum of $\mc G_{\sigma_{w,i}}$-invariant subcomplexes:
\begin{equation} \label{eq:9.9}
\text{gtr} \circ C_* (\tilde \sigma_\xi) : C_* \big( \mc F_\xi^{\mc G_{P,\delta}} 
\otimes_{e_{P,\delta} Z (\mc H \rtimes \Gamma )} e_{P,\delta} \mc H \rtimes \Gamma \big)
\to \bigoplus_{\sigma_{w,i} \prec (P,\delta)} 
C_* \big( \widehat{\mc O (T^{\sigma_{w,i}})}_{pr (\xi)} \big)^{\mc G_{\sigma_{w,i}}} .
\end{equation}
More precisely, the image of gtr $\circ C_* (\tilde \sigma_\xi)$ is the same as the formal
completion of 
\[
\bigoplus_{\sigma_{w,i} \prec (P,\delta)} \text{gtr} \circ C_* (\tilde \sigma )_{w,i} 
\big( C_* (\mc H \rtimes \Gamma) \big)
\]
at $\xi$. Since gtr $\circ C_* (\tilde \sigma)$ is a quasi-isomorphism 
(by Theorem \ref{thm:8.2}), and since 
\[
\mc H \rtimes \Gamma \mapsto \mc F_\xi^{\mc G_{P,\delta}} \otimes_{e_{P,\delta} 
Z (\mc H \rtimes \Gamma )} e_{P,\delta} \mc H \rtimes \Gamma 
\]
just has the effect of formally completing at pr$(\xi)$ and forgetting all analytic families 
with $\sigma_{w,i} \not\prec (P,\delta)$, \eqref{eq:9.9} is also a quasi-isomorphism. Thus
$\tilde \sigma_\xi$ induces
\[
HH_n \big( \mc F_\xi^{\mc G_{P,\delta}} \otimes_{e_{P,\delta} Z (\mc H \rtimes \Gamma )}
e_{P,\delta} \mc H \rtimes \Gamma \big) \cong \bigoplus_{\sigma_{w,i} \prec (P,\delta)} 
\widehat{\Omega^* (T^{\sigma_{w,i}})}_{pr (\xi)}^{\mc G_{\sigma_{w,i}}} .
\] 
By \eqref{eq:9.10}, \eqref{eq:9.8} and \eqref{eq:9.7} this means that, for all $\xi = 
(P,\delta,t) \in \Xi_{un} ,\: \tilde \sigma$ induces an isomorphism
\[
\mc F_\xi^{\mc G_{P,\delta}} \otimes_{e_{P,\delta} Z (\mc S \rtimes \Gamma )}
HH_n \big( e_{P,\delta} \mc S \rtimes \Gamma \big) \cong \bigoplus_{\sigma_{w,i} \prec 
(P,\delta)} \widehat{\Omega^* (T^{\sigma_{w,i}})}_{pr (\xi)}^{\mc G_{\sigma_{w,i}}} .
\]
In other words, the homomorphism of $Z (\mc S \rtimes \Gamma)$-modules \eqref{eq:9.5}
becomes an isomorphism upon formally completing with respect to any closed maximal ideal
of $Z (\mc S \rtimes \Gamma)$. Now the same argument as in the proof of Lemma \ref{lem:5.3}
shows that \eqref{eq:9.5} is itself an isomorphism. $\qquad \Box$
\vspace{2mm}

The following consequence of Theorem \ref{thm:9.1} can be proved in the same way as
Corollaries \ref{cor:7.2} and \ref{cor:8.3}.

\begin{cor}\label{cor:9.2}
The algebra homomorphism $\tilde \sigma_{un}$ induces isomorphisms
\[
\begin{array}{llll}
\hspace{-2mm} HC_* (\tilde \sigma_{un}) : & HC_* (\mc S \rtimes \Gamma) & \!\!\to\!\! &  
\bigoplus\limits_{w \in \langle W' \rangle} \bigoplus\limits_{i=1}^{c(w)} 
HC_* ( C^\infty (T_{un}^{\sigma_{w,i}} ))^{\mc G_{\sigma_{w,i}}} , \\
\hspace{-2mm} HP_{\text{ev/odd}} (\tilde \sigma) :  \hspace{-2mm} & 
HP_{\text{ev/odd}} (\mc S \rtimes \Gamma) & \!\!\to\!\! & 
\bigoplus\limits_{w \in \langle W' \rangle} \bigoplus\limits_{i=1}^{c(w)} H_{DR}^{\text{ev/odd}} 
(T_{un}^{\sigma_{w,i}} )^{\mc G_{\sigma_{w,i}}}  \cong  
H_{DR}^{\text{ev/odd}} (\widetilde T_{un})^{W'} .
\end{array} 
\]
\end{cor}

\subsection{Comparison of different parameters}

Let $k$ be real valued and let $q$ be positive. We will investigate what happens to the 
homology of the algebras 
$\mh H (\tilde{\mc R},k) \rtimes \Gamma ,\, \mc H (\mc R,q) \rtimes \Gamma$ and 
$\mc S (\mc R ,q) \rtimes \Gamma$ when we replace $k$ by $\ep k \; (\ep \in \C^\times)$ and 
$q$ by $q^\ep \; (\ep \in \C \setminus i \R)$. For the Schwartz algebras we must of 
course take $\ep$ real.

Recall the analytic families of representations $\pi_{\sigma_{w,i}}$ from Lemma \ref{lem:6.3} 
and $\pi_{\rho_w}$ from Corollary \ref{cor:6.5}. As discussed in the proof Theorem 
\ref{thm:6.6}, adjusting the parameter function by $\ep$ corresponds to applying 
$m_\ep^*$ to $\sigma_{w,i}$ and $\rho_w$. Let
\[
\begin{array}{lll}
\mh H (\tilde{\mc R},\ep k) \rtimes \Gamma & \to & 
\mc O (\mf t^{\rho_w}) \otimes \mr{End}_\C (V_{\rho_w}^\Gamma), \\
\mc H (\mc R,q^\ep) \rtimes \Gamma & \to & 
\mc O (T^{\sigma_{w,i}}) \otimes \mr{End}_\C (V_{\sigma_{w,i}}^\Gamma), \\
\mc S (\mc R,q) \rtimes \Gamma & \to & 
C^\infty (T_{un}^{\sigma_{w,i}}) \otimes \mr{End}_\C (V_{\sigma_{w,i}}^\Gamma)
\end{array}
\]
be the corresponding algebra homomorphisms. By Theorems \ref{thm:7.1}, \ref{thm:8.2} and 
\ref{thm:9.1} these induce isomorphisms
\begin{align*}
& HH_* (\mh H (\tilde{\mc R},\ep k) \rtimes \Gamma) \to \bigoplus_{w \in \langle W' \rangle} 
\Omega^* (\mf t^{\rho_w} )^{\mc G_{\rho_w}} \cong  \bigoplus_{w \in \langle W' \rangle} 
\Omega (\mf t^w)^{Z_{W'}(w)} , \\
&  HH_* (\mc H (\mc R,q^\ep) \rtimes \Gamma) \to \bigoplus_{w \in \langle W' \rangle} 
\bigoplus_{i=1}^{c(w)} \Omega^* (T^{\sigma_{w,i}} )^{\mc G_{\sigma_{w,i}}} \cong
\bigoplus_{w \in \langle W' \rangle} \Omega^* (T^w )^{Z_{W'}(w)} ,\\
& HH_* (\mc S (\mc R,q^\ep) \rtimes \Gamma) \to \bigoplus_{w \in \langle W' \rangle} 
\bigoplus_{i=1}^{c(w)} \Omega_{sm}^* (T_{un}^{\sigma_{w,i}})^{\mc G_{\sigma_{w,i}}} \cong
\bigoplus_{w \in \langle W' \rangle} \Omega_{sm}^* (T_{un}^w )^{Z_{W'}(w)} .
\end{align*}
The action of $Z (\mh H (\tilde{\mc R},\ep k) \rtimes \Gamma) \cong \mc O (\mf t)^{W'}$ on
$\Omega^* (\mf t^{\rho_w})$ is via the embedding
\[
\mf t^{\rho_w} \to \mf t : \lambda \mapsto cc(m_\ep^* (\rho_w)) + \lambda = 
\ep \, cc( \rho_w) + \lambda .                                  
\]
The translation factor $cc(m_\ep^* (\rho_w)) \in \mf t$ is a representative for the 
central character of $m_\ep^* (\rho_w)$. 
Similarly, the actions of $Z (\mc H (\mc R,q^\ep) \rtimes \Gamma) \cong \mc O (T )^{W'}$ on 
$\Omega^* (T^w)$ and on $\Omega^*_{sm} (T^w_{un})$ go via the embedding
\[
T^{\sigma_{w,i}} \to T : t \mapsto cc (m_\ep^* (\sigma_{w,i})) t . 
\]
Here $cc (m_\ep^* (\sigma_{w,i})) \in T_{un}$ is a representative for the central character 
of $\sigma_{w,i}$, and by \eqref{eq:6.weights} it lies in $Y \otimes S^1 q^{\ep \Z / 2}$.

Now we will generalize these results to the case $\ep = 0$. Recall the maps
\[
\zeta^\vee : G \big( \mh H (\tilde{\mc R},\ep k \big) \rtimes \Gamma \big) \to 
G \big( S (\mf t^*) \rtimes \Gamma \big)
\]
from Theorem \ref{thm:6.4} and
\[
\zeta^\vee : G \big( \mc  H (\mc R, q^\ep) \rtimes \Gamma \big) \to 
G \big( \C [X] \rtimes W' \big) 
\]
from Theorem \ref{thm:6.6}. These maps send analytic families to other analytic families, 
because they commute with parabolic induction.

\begin{thm} \label{thm:10.1}
\enuma{
\item Let $k$ be real valued and let $\ep \in \C^\times$. There exists a unique isomorphism 
\[
HH_* (\zeta^\vee) : HH_* (S (\mf t^*) \rtimes W') \to 
HH_* (\mh H (\tilde{\mc R},\ep k) \rtimes \Gamma)
\]
such that 
\[
HH_* (\pi_\rho) \circ HH_* (\zeta^\vee) = HH_* ( \zeta^\vee (\pi_\rho))
\]
for all analytic families of $\mh H (\tilde{\mc R}, \ep k) \rtimes \Gamma$-representations 
$\{ \pi_\rho (\lambda) \mid \lambda \in \mf t^\rho \}$.
\item Let $q$ be positive and let $\ep \in \C \setminus i\R$. There exists a unique isomorphism 
\[
HH_* (\zeta^\vee) : HH_* (\C [X] \rtimes W') \to HH_* (\mc H (\mc R,q^\ep) \rtimes \Gamma)
\]
such that 
\[
HH_* (\pi_\sigma) \circ HH_* (\zeta^\vee) = HH_* ( \zeta^\vee (\pi_\sigma))
\]
for all analytic families of $\mc H (\mc R,q^\ep) \rtimes \Gamma$-representations 
$\{ \pi_\sigma (t) \mid t \in T^\sigma \}$.
\item Let $q$ be positive. There exists a unique isomorphism 
\[
HH_* (\zeta^\vee) : HH_* (\mc S (X) \rtimes W') \to HH_* (\mc S (\mc R,q) \rtimes \Gamma)
\]
such that 
\[
HH_* (\pi_\sigma) \circ HH_* (\zeta^\vee) = HH_* ( \zeta^\vee (\pi_\sigma))
\]
for all analytic families of $\mc S (\mc R,q) \rtimes \Gamma$-representations 
$\{ \pi_\sigma (t) \mid t \in T_{un}^\sigma \}$.
} 
\noindent In all three settings the analogous statements for cyclic and periodic 
cyclic homology are also valid.
\end{thm}
\emph{Proof.}
We will only prove (a), the other parts follow in the same way. The construction of
the $\pi_{\rho_w}$ and the proof of Lemma \ref{lem:6.3} show that the analytic families
\begin{equation}\label{eq:10.1}
\{ \zeta^\vee (\pi_{\rho_w} (\lambda)) \mid w \in \langle W' \rangle , 
\lambda \in \mf t^{\rho_w} \}
\end{equation}
fulfill Corollary \ref{cor:6.5} for $\mh H (\tilde{\mc R},0) \rtimes \Gamma = S (\mf t^*) 
\rtimes W'$. We remark that the representations \eqref{eq:10.1} can be reducible, even
for generic $\lambda$. But since this does not affect the arguments in Section 
\ref{sec:HGHA}, we may ignore it. In particular Theorem \ref{thm:7.1} holds for
$S(\mf t^*) \rtimes W'$ with \eqref{eq:10.1}. Now the condition of the theorem enforces
that $HH_* (\zeta^\vee)$ is the composite isomorphism 
\[
HH_* (S(\mf t^*) \rtimes W') \longrightarrow \bigoplus_{w \in \langle W' \rangle} \Omega^* 
(\mf t^{\rho_w} )^{\mc W'_{\rho_w}} \xrightarrow{\: HH_* (\tilde \rho )^{-1} \:} 
HH_* (\mh H (\tilde{\mc R},\ep k) \rtimes \Gamma ) .
\]
For all $w \in \langle W' \rangle$ this map satisfies
\[
HH_* (\pi_{\rho_w}) \circ HH_* (\zeta^\vee) = HH_* ( \zeta^\vee (\pi_{\rho_w})) .
\]
Thus it remains to check compatibility with an arbitrary analytic family of
$\mh H (\tilde{\mc R},\ep k) \rtimes \Gamma$-representations $\{ \pi_\rho (\lambda) \mid
\lambda \in \mf t^{\rho} \}$. Let $\lambda$ be a generic point of $\mf t^\rho$.
By Corollary \ref{cor:6.5}.c we can find $n_\rho \in \N, n_w \in \Z$ and 
$\mu_{w,\lambda} \in \mf t^{\rho_w}$ such that
\begin{equation} \label{eq:10.2}
n_\rho \pi_\rho (\lambda) = \sum_{w \in \langle W' \rangle} n_w \pi_{\rho_w}(\mu_{w,\lambda})
\quad \text{in} \quad G \big( \mh H (\tilde{\mc R},\ep k) \rtimes \Gamma \big) .
\end{equation}
The genericity of $\lambda$ implies $n_w = 0$ unless $W' \mf t^{\rho_w} \supset W' \mf t^\rho$.
Comparing the central characters also shows that there exist $\tilde w \in W'$ with
\[
\tilde \lambda + cc (\pi_\rho (0)) = \mu_{w,\lambda} + cc (\pi_{\rho_w}(0)) ,
\]
for suitable representatives $cc (\pi)$ of the central character of the representations
$\pi$ under consideration. Hence we may take
\[
\mu_{w,\lambda} = \tilde w \lambda + cc(\pi_\rho (0)) - cc (\pi_{\rho_w} (0)) =: 
\tilde \lambda + \mu_w .
\]
Since $\lambda$ was generic, with this choice of $\mu_{w,\lambda}$ \eqref{eq:10.2} becomes
valid for all generic $\lambda \in \mf t\rho$, which then by continuity extends to the 
whole of $\mf t^\rho$. Next we bring the terms $n_w < 0$ to the other side and obtain
\[
n_\rho \pi_\rho (\lambda) + \sum_{w \in \langle W' \rangle_-} -n_w \pi_{\rho_w}(\tilde \lambda 
+ \mu_w) = \sum_{w \in \langle W' \rangle_+} n_w \pi_{\rho_w}(\tilde \lambda + \mu_w) .
\]
This is an equality in $G (\mh H (\tilde{\mc R},\ep k) \rtimes \Gamma)$, but the 
representations on the both sides need not be isomorphic. Since Hochschild homology does not 
distinguish direct sums from nontrivial extensions \cite[Theorem 1.2.15]{Lod}, we may 
nevertheless treat it as an equivalence of 
$\mh H (\tilde{\mc R},k) \rtimes \Gamma$-representations. Consequently
\[
n_\rho HH_* (\pi_\rho ) + \sum_{w \in \langle W' \rangle_-} -n_w 
HH_* (\pi_{\rho_w})\big|_{\mu_w + \tilde w \mf t^\rho} = 
\sum_{w \in \langle W' \rangle_+} n_w HH_* (\pi_{\rho_w}) \big|_{\mu_w + \tilde w \mf t^\rho}
\]
as maps 
\[
HH_* \big( \mh H (\tilde{\mc R},\ep k) \rtimes \Gamma \big) \to HH_* (\mc O (\mf t^\rho)) 
\cong \Omega^* (\mf t^\rho) .
\]
Now he aforementioned compatibility of $HH_* (\zeta^\vee)$ with the $\pi_{\rho_w}$ and the
additivity of $\zeta^\vee$ show that 
\[
HH_* (\pi_\rho) \circ HH_* (\zeta^\vee) = HH_* ( \zeta^\vee (\pi_\rho)) ,
\]
as required. The analogous statements for (periodic) cyclic homology can be proved in the 
same way, using Corollary \ref{cor:7.2}. $\qquad \Box$
\vspace{2mm}

In \cite[Theorem 4.4.2]{Sol-Irr} the author constructed a homomorphism of Fr\'echet algebras
\[
\zeta_0 : \mc S (X) \rtimes W' \to \mc S (\mc R,q) \rtimes \Gamma .
\]
It is related to Theorem \ref{thm:6.2} by \cite[Corollary 4.4.3]{Sol-Irr}: $\zeta^\vee (\pi)
\cong \pi \circ \zeta_0$ for all irreducible $\mc S (\mc R,q) \rtimes \Gamma$-representations. 
Furthermore it was shown in \cite[Theorems 5.1.4 and 5.2.1]{Sol-Irr} that $\zeta_0$ induces 
isomorphisms on topological $K$-theory and on periodic cyclic homology.
With Theorem \ref{thm:10.1} we can improve on this:

\begin{prop} \label{prop:10.2}
Let $q$ be positive.
\enuma{
\item $HH_* (\zeta_0) = HH_* (\zeta^\vee) : HH_* (\mc S (X) \rtimes W') \to
HH_* (\mc S (\mc R,q) \rtimes \Gamma )$.
\item There exists a natural commutative diagram
\[
\begin{array}{ccc}
HH_* (\C [X] \rtimes W') & \to & HH_* (\mc S (X) \rtimes W') \\
\downarrow \scs{HH_* (\zeta^\vee)} & & \downarrow \scs{HH_* (\zeta_0)} \\
HH_* (\mc H (\mc R,q) \rtimes \Gamma) & \to & HH_* (\mc S (\mc R,q) \rtimes \Gamma)
\end{array}
\]
in which the horizontal maps are induced by inclusions of algebras and the
vertical maps are isomorphisms.
\item The analogues of (a) and (b) for cyclic and periodic cyclic homology are
also valid.
}
\end{prop}
\emph{Proof.}
By \cite[Lemma 4.2.3 and Theorem 4.4.2.e]{Sol-Irr}
\[
\pi_\sigma (t) \circ \zeta_0 = \zeta^\vee (\pi_\sigma (t))
\]
for every tempered analytic family of $\mc H (\mc R,q) \rtimes \Gamma$-representations
$\{ \pi_\sigma (t) \mid t \in T^\sigma_{un} \}$. Thus (a) follows from the unicity part
of Theorem \ref{thm:10.1}.c. That and Theorem \ref{thm:10.1} directly imply (b) and (c).
$\qquad \Box$ 

\end{document}